\documentclass[11pt,leqno,oneside, letterpaper]{amsart}

\usepackage[margin=1.15in]{geometry}
\usepackage{amsmath,amsfonts,amssymb,amscd,amsthm,amsbsy,epsf,graphics} 
\usepackage{stmaryrd}
\usepackage{setspace}
\usepackage{comment}
\usepackage{color}
\usepackage[dvipsnames]{xcolor}
\usepackage[colorlinks,linkcolor=blue,citecolor=blue]{hyperref} 
\usepackage{epsfig} 
\usepackage{tikz-cd}
\usepackage{bm}
\usepackage{graphicx}  
\usepackage[caption=false]{subfig}
\usepackage{lipsum}
\usepackage{enumitem} 

\usepackage[square,numbers]{natbib}

\makeatletter
\def\l@subsection{\@tocline{2}{0pt}{4pc}{5pc}{}}
\makeatother
 
\let\oldtocsection=\tocsection
\let\oldtocsubsection=\tocsubsection
\let\oldtocsubsubsection=\tocsubsubsection
 
\renewcommand{\tocsection}[2]{\hspace{0em}\oldtocsection{#1}{#2}}
\renewcommand{\tocsubsection}[2]{\hspace{0em}\oldtocsubsection{#1}{#2}}
\renewcommand{\tocsubsubsection}[2]{\hspace{2em}\oldtocsubsubsection{#1}{#2}}

\newtheorem{thm}{Theorem}[section]
\newtheorem*{thm*}{Theorem}
\newtheorem{cor}[thm]{Corollary}
\newtheorem{lemma}[thm]{Lemma}
\newtheorem{prop}[thm]{Proposition}

\newtheorem{innercustomthm}{Theorem}
\newenvironment{customthm}[1]
  {\renewcommand\theinnercustomthm{#1}\innercustomthm}
  {\endinnercustomthm}

\newtheorem{innercustomcor}{Corollary}
\newenvironment{customcor}[1]
{\renewcommand\theinnercustomcor{#1}\innercustomcor}
{\endinnercustomcor}

\theoremstyle{definition}

\newtheorem{definition}[thm]{Definition} 
\newtheorem{example}[thm]{Example}

\newtheorem{remark}[thm]{Remark} 

\newtheorem*{notation}{Notation}
\newtheorem*{organization}{Plan of the paper}
\newtheorem*{acknowledgement*}{Acknowledgements}

\newtheoremstyle{cases}
  {12pt plus 6 pt}
  {2pt}
  {\bfseries}   
  {}
  {\bfseries}
  {.}
  {.5em}
  {}
\theoremstyle{cases}

\numberwithin{subcase}{case} 
\numberwithin{subsubcase}{subcase}
\numberwithin{equation}{subsection} 

\newcommand{\G}{\text{Homeo}_+(S^1)}
\newcommand{\tG}{\text{Homeo}_{\mathbb{Z}}(\mathbb{R})}

\graphicspath{{figures/}}

\title{Recalibrating $\mathbb{R}$-order trees and $\mbox{Homeo}_+(S^1)$-representations of link groups}

\author[Steven Boyer]{Steven Boyer} 
\thanks{Steven Boyer was partially supported by NSERC grant RGPIN 9446-2008}
\address{D\'epartement de Math\'ematiques, Universit\'e du Qu\'ebec \`a Montr\'eal, 201 President Kennedy Avenue, Montr\'eal, Qc., Canada H2X 3Y7.}
\email{boyer.steven@uqam.ca}
\urladdr{http://www.cirget.uqam.ca/boyer/boyer.html}

\author[Cameron McA. Gordon]{Cameron McA. Gordon}
\address{Department of Mathematics, University of Texas at Austin, 1 University Station, Austin, TX 78712, USA.}
\email{gordon@math.utexas.edu}

\author[Ying Hu]{Ying Hu}
\address{Department of Mathematical and Statistical Sciences, University of Nebraska Omaha, 6001 Dodge Street, Omaha, NE 68182-0243, USA.}
\email{yinghu@unomaha.edu}
\urladdr{https://yinghu-math.github.io}


\date{}

\begin{document}
\begin{abstract}
In this paper we study the left-orderability of $3$-manifold groups using an enhancement, called recalibration, of Calegari and Dunfield's ``flipping'' construction, used for modifying $\G$-representations of the fundamental groups of closed $3$-manifolds. The added flexibility accorded by recalibration allows us to produce $\G$-representations of hyperbolic link exteriors so that a chosen element in the peripheral subgroup is sent to any given rational rotation. We apply these representations to show that the branched covers of families of links associated to epimorphisms of the link group onto a finite cyclic group are left-orderable. This applies, for instance, to fibered hyperbolic strongly quasipositive links. Our result on the orderability of branched covers implies that the degeneracy locus of any pseudo-Anosov flow on an alternating knot complement must be meridional, which generalizes the known result that the fractional Dehn twist coefficient of any hyperbolic fibered alternating knot is zero. Applications of these representations to order-detection of slopes are also discussed in the paper.
\end{abstract}

\maketitle

\section{Introduction} 
\label{sec: intro}
The study of left-orderability of $3$-manifold groups was initiated in \cite{BRW05} and has proven to be closely tied to various aspects of the topology of $3$-manifolds \cite{BGW13,CR16}. By \cite[Theorem 1.1]{BRW05}, the fundamental group of an irreducible $3$-manifold is left-orderable if it has a nontrivial ${\text{Homeo}}_+(\mathbb{R})$-representation. However, it is difficult to construct ${\text{Homeo}}_+(\mathbb{R})$-representations directly. Instead, most results on the left-orderability of $3$-manifold groups are consequences of the existence of $\G$-representations of a link group with vanishing Euler class, which implies that the action on the circle can be lifted to an action on the real line. For instance, see \cite{Hu15, Tran15, Gor17, CD18, BH19, Dun20, Turner21, DuRa21, BGH21}.

Various structures on a closed $3$-manifold $W$ such as foliations, laminations, and flows often result in interesting representations $\rho: \pi_1(W) \to \text{Homeo}_+(S^1)$. The $\text{Homeo}_+(S^1)$-representations associated to co-oriented taut foliations (\cite{Thurston}; see \cite[Theorem 6.2]{CD03}) were used in \cite{BH19,BGH21} to study the left-orderablity of $3$-manifold groups, as have the $\G$-representations arising from very full laminations with orderable cataclysms (\cite[Theorem 3.8]{CD03}) and from pseudo-Anosov flows (\cite{Fen12}, \cite{CD03}). 

We refer the reader to \S \ref{subsec: psan flows} and \S \ref{subsec: degeneracy loci} for the definition of pseudo-Anosov flows, stable and unstable laminations as well as other related concepts. A good example of a pseudo-Anosov flow to keep in mind is the suspension of a pseudo-Anosov homeomorphism $\varphi$ of a surface of finite type. In this example, the stable and unstable laminations are the suspensions of the stable and unstable invariant laminations of $\varphi$, and both are very full.

In this article, we continue this theme of studying the left-orderability of $3$-manifold groups by utilising circle actions from pseudo-Anosov flows and very full laminations, but with two new ingredients: 
\begin{enumerate}
    \item We extend to an orbifold setting the constructions of Calegari-Dunfield and Fenley of $\G$-representations associated to pseudo-Anosov flows and to very full laminations. (See Proposition \ref{prop: pa implies co orbifold} and Theorem \ref{thm: pa implies co orbifold 2}).  
   \item We reformulate and generalise the flipping operation described in \cite[\S 3.3]{CD03} as a purely combinatorial operation on cyclically-ordered $\mathbb{R}$-order trees that we call {\it recalibration}. This allows for much more flexibility in producing $\G$-representations (see \S \ref{subsec: pa flow action on the tree}). 
\end{enumerate}
Combining these ideas allows us to construct many nontrivial representations of link groups with prescribed behaviour on the peripheral subgroups. To describe this more precisely, we first introduce some notation. 

\begin{notation}
Let $L = K_1 \cup \cdots \cup K_m$ be a link in a closed, orientable 3-manifold $W$ with closed tubular neighbourhood $N(L)$. The {\it complement} of $L$ in $W$ is the open manifold $C(L) = W \setminus L$ and the {\it exterior} of $L$ is the compact manifold $X(L) = W \setminus \mbox{int}(N(L))$. We often identify $C(L)$ with the interior of $X(L)$. We use $T_i = \partial N(K_i) \subset \partial X(L)$ to denote the $i^{th}$ boundary component of $X(L)$. An essential oriented simple closed curve $\alpha$ on $T_i$ determines a slope on $T_i$, by forgetting its orientation, and a class in $\pi_1(T_i) = H_1(T_i)$. We also use $\alpha_i$ to denote the corresponding class in $\pi_1(X(L))$, which is well-defined up to conjugation.  For each $i$, we use $\mu_i$ to denote an oriented meridional curve of $N(K_i)$. Links are unoriented unless otherwise stated or clear from the context. If $L$ and $W$ are both oriented, then the $\mu_i$ are oriented positively with respect to the orientations on $L$ and $W$.
\end{notation}

In the theorem below, $\delta_i(\Phi_0)$ is the degeneracy locus on $T_i$ of a pseudo-Anosov flow $\Phi_0$ on $C(L)$. For instance, if $L$ is fibred with pseudo-Anosov monodromy and $\Phi_0$ is the suspension flow of the pseudo-Anosov representative of the monodromy restricted to $C(L)$, thought of as the interior of $X(L)$,  then $\delta_i(\Phi_0)$ is a collection of parallel simple closed curves on $T_i$ consisting of half of the closed flow lines on $T_i$. See \S \ref{subsec: degeneracy loci} for more details. 

\begin{customthm}{\ref{thm: reps result intro}}
Let $L = K_1 \cup \cdots \cup K_m$ be a link in an orientable $3$-manifold $W$ whose complement admits a pseudo-Anosov flow $\Phi_0$. For each $i$, fix an oriented essential simple closed curve $\alpha_i$ on $T_i$ and an integer $n_i \geq 1$ so that $n_i|\alpha_i \cdot \delta_i(\Phi_0)| \geq 2$. Then for any integer $a_i$ coprime with $n_i$, there is a homomorphism $\rho: \pi_1(X(L)) \to \mbox{Homeo}_+(S^1)$ with non-cyclic image such that $\rho(\alpha_i)$ is conjugate to rotation by $2\pi a_i/n_i$ for each $i$.
\end{customthm}

The theorem shows that given any multislope $(\alpha_1, \ldots, \alpha_m)$ on $\partial X(L)$ such that $\delta_i(\Phi_0)$ is not a multiple of $\alpha_i$ and any $m$-tuple $(u_1, u_2, \ldots, u_m)$ of non-trivial roots of unity, there is a  homomorphism $\rho: \pi_1(X(L)) \to \mbox{Homeo}_+(S^1)$ with non-cyclic image such that for each $i$, $\rho(\alpha_i)$ is conjugate to the rotation determined by $u_i$. 

\subsection*{Applications to order-detection}
The notion of slope detection was first introduced by Boyer and Clay in \cite{BC} to prove the $L$-space conjecture for graph manifolds, and further used in \cite{BGH21} by the authors of the present article to study the $L$-space conjecture for general toroidal manifolds (\S \ref{subsec: app slope detection}). One of the main results of the latter paper was to show that meridional slopes (i.e. slopes that are distance one from the longitudinal slope) on the boundaries of integer homology solid tori are order-detected, a result with many consequences. For instance, the following results from \cite{BGH21} and \cite{BGH2} are obtained as applications of this fact: 
\begin{enumerate}
    \item The fundamental groups of cyclic branched covers of prime satellite links are all left-orderable \cite{BGH21}, which settles the left-orderability version of the Gordon-Lidman Conjecture \cite{GL14}.
    \item The fundamental groups of toroidal integer homology spheres are left-orderable \cite{BGH21}, which is the left-orderable analog of a result of \cite{Eft18,HRW1}  on $L$-spaces. 
    \item The fundamental group of any non-meridional Dehn filling on a composite knot is left-orderable \cite{BGH2}, which is predicted by the $L$-space conjecture and the fact that $L$-space knots are prime \cite{Krc15}.
\end{enumerate}

The proof that meridional slopes are order-detected in \cite[Theorem 1.3(2), Corollary 1.4]{BGH21} involves a detailed analysis of the dynamics of Thurston's universal circle action associated with taut foliations, and is long and technical. 
In \S \ref{subsec: app slope detection}, we give a very quick proof that meridians of certain hyperbolic knots are order-detected (see Theorem \ref{thm: meridian detection}) using the representations constructed in Theorem \ref{thm: reps result intro}. 

\begin{customthm}{\ref{thm: meridian detection}}
Let $K$ be a hyperbolic knot in an integer homology sphere $W$, such that the complement of $K$ admits a pseudo-Anosov flow whose degeneracy locus is meridional. Then the meridional slope of $K$ is order-detected. 
\end{customthm}

Theorem \ref{thm: meridian detection} demonstrates the potential of these representations in studying order-detected slopes. We plan to undertake a fuller investigation of this in a future paper.

\subsection*{Applications to left-orderability of branched covers}
We consider a general class of cyclic branched covers defined as follows.

\begin{definition}
Let $L = K_1 \cup ... \cup K_m$ be a link in an integer homology sphere $W$ and $\psi : \pi_1(X(L)) \to \mathbb Z/n$ an epimorphism such that $\psi(\mu_i) \ne 0$ for each $i$. We denote by $\Sigma_\psi(L) \to W$ the $n$-fold cyclic cover $\Sigma_\psi(L) \to W$ branched over $L$ associated to the $n$-fold cyclic cover $X_\psi(L) \to X(L)$ determined by $\psi$. 
\end{definition}

Following standard notational conventions, if $L$ is oriented we use $\Sigma_n(L)$ to denote the {\it canonical cyclic branched cover} of $L$ associated with the epimorphism $\pi_1(X(L)) \to \mathbb{Z}/n$ that sends each correspondingly oriented $\mu_i$ to $1$ (mod $n$). 

Motivated by a result of Ba and Clay, \cite[Theorem 2.6]{BaC}, our next theorem shows how a representation of the sort given by Theorem \ref{thm: reps result intro} can be used to deduce the left-orderability of the fundamental group of an associated cyclic branched cover $\Sigma_\psi(L)$. 
This result greatly generalizes \cite[Theorem 3.1]{Hu15} and \cite[Theorem 6]{BGW13}, which have been widely used for proving the left-orderability of the fundamental group of cyclic branched covers of knots. 

\begin{customthm}{\ref{thm: locbc intro}}
Let $L = K_1 \cup \cdots \cup K_m$ be a prime link in an integer homology $3$-sphere whose exterior is irreducible. Suppose that $\rho: \pi_1(X(L)) \to \mbox{Homeo}_+(S^1)$ is a representation with non-cyclic image such that $\rho(\mu_i)$ is conjugate to rotation by $2\pi a_i/n$ for some $a_i, n\in \mathbb{Z}$, where $n \geq 2$. If the induced homomorphism $\psi: \pi_1(X(L)) \to \mathbb{Z}/n$ which sends $\mu_i$ to $a_i$ $(${\rm{mod}} $n)$ is an epimorphism, then $\pi_1(\Sigma_\psi(L))$ is left-orderable. 
\end{customthm}

Combining Theorem \ref{thm: reps result intro} with Theorem \ref{thm: locbc intro}, we found, somewhat surprisingly, many examples of hyperbolic links, for which the fundamental  group of any cyclic branched cover is left-orderable (see Theorem \ref{thm: lo cbcs intro} below), even though for a given covering index these branched covers can be very different as topological spaces. Since links whose  exteriors admit pseudo-Anosov flows are prime (cf. Remark \ref{rem: primeness of links}), Theorem \ref{thm: lo cbcs intro} is a corollary of these results; the proof is given in \S \ref{subsec: lo and link orientation}.

\begin{thm} 
\label{thm: lo cbcs intro}
Let $L$ be a link in an integer homology sphere whose complement admits a pseudo-Anosov flow none of whose degeneracy loci are meridional. Then the fundamental group of any $n$-fold cyclic branched cover $\Sigma_{\psi}(L)$ of $L, n \ge 2$, is left-orderable. 
\end{thm}

\begin{example}
In \S \ref{subsec; egs}, we give many examples of links which satisfy the hypothesis of Theorem \ref{thm: lo cbcs intro}. These include: 
\begin{enumerate}
\setlength\itemsep{0.3em}
\item Hyperbolic links that can be oriented to be fibered with nonzero fractional Dehn twist coefficient (Theorem \ref{thm: hyperbolic fibre sigman $LO$ intro}). Most interestingly, links that can be oriented to be fibered and strongly quasipositive belong to this family (Corollary \ref{cor: hyperbolic fibred sqp links intro}). 
\item Links that can be oriented to be the closures of certain pseudo-Anosov braids (see Theorem \ref{thm: closed braids}, Corollary \ref{cor: fdtc braid}, Proposition \ref{prop: example of braids}, and Theorem \ref{thm: example of braids}).
\end{enumerate}
\end{example}

The following corollary is a special case of the first of these two examples.

\begin{customcor}{\ref{cor: l-space knots branched covers}}
If $K$ is an $L$-space knot then $\pi_1(\Sigma_n(K))$ is left-orderable for all $n \ge 2$ if and only if $K$ is not $T(3, 4), T(3, 5)$, or $T(2, 2q+1)$ for some $q \geq 1$.
\end{customcor}

Using \cite{BBG19-1} and \cite{FRW22} we obtain the analogous statement with ``$\pi_1(\Sigma_n(K))$ is left-orderable" replaced by ``$\Sigma_n(K)$ is not an $L$-space", provided $n \ge 3$. This leaves open the interesting question, due to Allison Moore, asking whether the double branched cover of a hyperbolic $L$-space knot can ever be an $L$-space. 

We use $L^\mathfrak{0}$ to denote $L$ endowed with an orientation $\mathfrak{o}$. Varying the values of $\psi(\mu_i)$ in Theorem \ref{thm: lo cbcs intro} over all possible choices of $\pm 1$ (mod $n$) yields the following corollary. 

\begin{cor} 
\label{cor: $LO$ branched cover with all orientations intro}
Let $L = K_1 \cup \cdots \cup K_m$ be a link in an integer homology $3$-sphere $W$ whose complement admits a pseudo-Anosov flow none of whose degeneracy loci are meridional. Then $\pi_1(\Sigma_n(L^\mathfrak{o}))$ is left-orderable for all $n \ge 2$ and all orientations $\mathfrak{o}$ on $L$,
\end{cor}

To put Corollary \ref{cor: $LO$ branched cover with all orientations intro} in context, we remark that known results suggest that if the double branched cover of an oriented link in an integer homology $3$-sphere has a left-orderable fundamental group, then so does $\Sigma_n(L)$ for all $n \geq 2$. On the other hand, $\Sigma_2(L)$ is independent of the orientation on $L$, so we expect that 
the left-orderability of $\pi_1(\Sigma_2(L))$ implies that of $\pi_1(\Sigma_n(L^\mathfrak{o}))$ for any orientation $\mathfrak{o}$ on $L$ and $n \geq 2$. To the best of our knowledge, Corollary \ref{cor: $LO$ branched cover with all orientations intro} is the first general result confirming this type of behaviour for hyperbolic links. See \S \ref{subsec: lo and link orientation} for a more detailed discussion.

We complement Corollary \ref{cor: $LO$ branched cover with all orientations intro} by showing that there are pairs of oriented links $L, L'$, which are equivalent as unoriented links, such that $\pi_1(\Sigma_n(L))$ is non-left-orderable for all $n \ge 2$, while $\pi_1(\Sigma_n(L'))$ is left-orderable for all $n \ge 3$. See \S \ref{subsec: egs}.

\subsection*{Applications to degeneracy loci}
Lastly, we present a surprising application of our left-orderability results to degeneracy loci of pseudo-Anosov flows. Independent work of Gabai and Mosher (\cite{Mosher96}; also see \cite{LT}) showed that pseudo-Anosov flows exist in the complement of any hyperbolic link in a closed, orientable $3$-manifold. On the other hand, the degeneracy loci of these flows are difficult to compute in general, although they are known in certain situations. For instance,  \cite[Corollary 1.7]{BNS22} combines with \cite[Theorem 1.2]{Hed10} and \cite[Corollary 7.3]{BBG19-1} to show that the degeneracy locus of the monodromy of a fibred hyperbolic alternating knot is always meridional. We significantly extend this fact by showing : 

\begin{customcor}{\ref{cor: degeneracy loci hyp alt knots}}
The degeneracy locus of any pseudo-Anosov flow on the complement of an alternating knot is meridional.  
\end{customcor}

\begin{organization}
Section \ref{sec: preliminaries} outlines background material on Fenley's asymptotic circle representations from pseudo-Anosov flows on closed $3$-manifolds, essential laminations and degeneracy loci, and rotation and translation numbers. In \S \ref{sec: univ circs and locbcs hyp} we extend Fenley's asymptotic circle representation to the fundamental groups of closed orientable $3$-orbifolds with cyclic isotropy (Proposition \ref{prop: pa implies co orbifold}). These representations are used to deal with a degenerate case in our arguments in \S \ref{sec: psA flows and $LO$ cbcs}. In \S \ref{sec: psA flows and $LO$ cbcs} we define the recalibration operation on cyclically ordered $\mathbb{R}$-order trees, and use its flexibility to produce $\mbox{Homeo}_+(S^1)$-valued representations of link groups with prescribed behaviour on the peripheral subgroups (Theorem \ref{thm: reps result intro}). In \S \ref{sec: reps and locbcs}, we consider the Euler classes of these representations and prove Theorem \ref{thm: locbc intro}. The applications discussed above are proved in \S \ref{sec: examples}. Finally, in \S \ref{sec: applications}, we discuss our results in the context of the L-space conjecture. 
\end{organization}
 
\begin{acknowledgement*}
The authors would like to thank John Baldwin, who contributed to initial discussions of the material in this paper, Sergio Fenley, for providing background information on his asymptotic circle, and Adam Clay, for telling us of his work with Idrissa Ba on co-cyclic left-orderable subgroups of circularly-ordered groups. 
\end{acknowledgement*}

\section{Preliminaries}  
\label{sec: preliminaries}
In \S \ref{subsec: psan flows} and \S \ref{subsec: degeneracy loci}, we cover some basic concepts related to pseudo-Anosov flows. We refer the readers to \cite[\S 6.6]{Cal07} for a more detailed account of the material. We explain the relationship between degeneracy loci and fractional Dehn twist coefficients in \S \ref{subsec: fdtc}.  In \S \ref{subsec: fenley's uc}, we survey Fenley's results on the existence of a circle action given a pseudo-Anosov flow on a closed $3$-manifold. Finally, in \S \ref{subsec: rot and trans numbers}, we briefly review the rotation and translation numbers of elements in $\G$ and $\tG$.

\subsection{Pseudo-Anosov flows}
\label{subsec: psan flows} 
An {\it Anosov flow} $\Phi_t: M\times \mathbb{R} \rightarrow M$ on a $3$-manifold $M$ is a flow which preserves a continuous splitting of the
tangent bundle $TM = E^s\oplus \frac{\partial}{\partial t} \oplus E^u $. Moreover, there are constants $\mu_0 \geq 1$ and $\mu_1 > 0$ so that
$$\Vert d\Phi_t(v) \Vert \leq \mu_0 e^{- \mu_1 t} \Vert v\Vert,  $$ 
$$ \Vert d\Phi_{-t}(w) \Vert \leq \mu_0 e^{- \mu_1 t} \Vert w\Vert, $$
for any $v\in E^s$, $w\in E^u$ and $t\geq 0$. So an Anosov flow contracts vectors along $E^s$ and expands vectors along $E^u$. By definition, $d\Phi_t$ is a hyperbolic map on $E^s\oplus E^u$. 

A flow $\Phi_t$ on a $3$-manifold $M$ is {\it pseudo-Anosov} if  it is Anosov away from a finite number of {\it pseudo-hyperbolic}  periodic orbits. That is, $d\Phi_t$ restricted to the normal bundle of the flow has a {\it pseudo-hyperbolic} singularity at each point along these singular orbits. The archetypical example of a pseudo-Anosov flow is the suspension of a pseudo-Anosov homeomorphism of a connected, orientable surface, in which the periodic orbits obtained from the singular points of the homeomorphism's invariant singular foliations are pseudo-hyperbolic. 

In this article, we include Anosov flows when we refer to pseudo-Anosov flows, and will consider pseudo-Anosov flows on closed manifolds as well as link complements in a closed manifold, though in the latter case we restrict the behaviour of the flow in the ends of the link complement as follows. 

Let $C(L) = W\setminus L$ be the complement of a link $L$ in a closed, connected, orientable $3$-manifold $W$. For a flow on $C(L)$ to be pseudo-Anosov, in addition to the definition above, we also require the dynamics of each end of $C(L)$ to be that of a neighbourhood of a pseudo-hyperbolic orbit in a pseudo-Anosov flow with the orbit removed. The archetypical example of such a flow is the suspension of a pseudo-Anosov homeomorphism of a cusped surface of finite type. 

To simplify notation, we will drop the subscript $t$ from the notation of a flow from now on.

\subsection{Essential laminations and degeneracy loci}
\label{subsec: degeneracy loci}
Let $\Phi$ be a pseudo-Anosov flow on a $3$-manifold $M$.
By the stable manifold theorem \cite{HPS77}, both $E^s\oplus \frac{\partial}{\partial t}$ and $E^u\oplus \frac{\partial}{\partial t}$ are integrable, which results in two singular foliations on $M$ that are invariant under the flow. They are called the weak stable foliation and the weak unstable foliation respectively. 

A lamination on a $3$-manifold $M$ is a foliation on a closed subset of $M$ by surfaces. So given a pseudo-Anosov flow on $M$, by blowing air into the singular leaves of its weak stable and unstable invariant foliations (i.e. replacing the leaves by small regular neighbourhoods whose interiors are removed), one obtains two invariant {\it essential} laminations on $M$, which are called the stable and unstable laminations of the pseudo-Anosov flow. Next, we give a description of the topological structure of the invariant laminations of a pseudo-Anosov flow which is sufficient for our purposes. For more general background material on essential laminations, we refer the reader to \cite{GO89}.

Let $\Lambda$ be a lamination on $M$. A {\it complementary region} of $\Lambda$ is a component of the completion of $M\setminus \Lambda$ under a path metric on $M$. We call a lamination on a $3$-manifold $M$ {\it very full} (\cite[Definition 6.42]{Cal07}) if each complementary region is homeomorphic to either an ideal polygon bundle over $S^1$ or a once punctured ideal polygon bundle over $S^1$. In the first case, we require the ideal polygon to have at least $2$ ideal vertices. We follow the terminology in \cite{Mosher96} and call a complementary region of the first type a {\it pared solid torus} and of the second type {\it a pared torus shell}. By construction, the invariant laminations of a pseudo-Anosov flow on $M$ are very full.

Let $\Phi_0$ be a pseudo-Anosov flow on the complement $C(L)$ of a link $L = K_1 \cup K_2 \cup \cdots \cup K_m$ in $W$ and $\Lambda$ be the stable lamination of $\Phi_0$. Then for each $i$, there is a pared torus shell complementary region of $\Lambda$, corresponding to the missing $i^{th}$ component $K_i$ of the link $L$. Note that the pared torus shell is homeomorphic to $N(K_i)\setminus K_i$ with a collection of parallel simple closed curves on $T_i = \partial N(K_i)$ removed. The isotopy class of this collection of simple closed curves is called the {\it degeneracy locus} of the pseudo-Anosov flow on $T_i$ and is denoted by $\delta_{T_i}(\Phi_0)$ or, more simply, by $\delta_i(\Phi_0)$. This notion was first introduced in \cite{GO89} for the suspension flow of a pseudo-Anosov homeomorphism, where it was called the degenerate curve. These degeneracy loci are important in determining if an essential lamination remains essential after Dehn filling (\cite[Theorem 5.3]{GO89}).

We say that the degeneracy locus is {\em meridional} on the $i^{th}$ component $T_i$ if $\delta_i(\Phi_0)$ is a union of meridional curves of $N(K_i)$. The slope determined by $\delta_i(\Phi)$ on $T_i$ is the slope determined by a connected component of $\delta_i(\Phi)$.

\subsection{Degeneracy loci and fractional Dehn twist coefficients}
\label{subsec: fdtc}
Let $L = K_1\cup \cdots \cup K_m$ be a hyperbolic fibred link in an oriented $3$-manifold $W$ with monodromy $h$ and fibre $F$ (so $(F,h)$ is an open book decomposition of $W$). Since $L$ is hyperbolic,  $h$ is freely isotopic to a pseudo-Anosov homeomorphism $\varphi$ of $F$ \cite{Thurston98}. The suspension flow $\Phi$ of $\varphi$ on the link exterior $X(L)$ restricts to a  pseudo-Anosov flow on the link complement $int(X(L)) \cong C(L)$ which we denote by $\Phi_0$. In this case, the degeneracy loci of $\Phi_0$ can be described more precisely, as we explain next. 

It follows from the properties of pseudo-Anosov homeomorphisms that the flow $\Phi$ has an even number of periodic orbits on each boundary component of $X(L)$, half from the repelling periodic points of $\varphi$ on $\partial F$ and half from the attracting periodic points. Then on each boundary component, the degeneracy locus of $\Phi_0$ is the union of half of the periodic orbits of $\Phi$ on $\partial X(L)$. 

For each $i = 1, \cdots, m$, let $T_i$ denote the boundary component $\partial N(K_i)$ of $X(L)$, $\mu_i$ the meridional class given by $h$, and $\lambda_i$ the slope on $T_i$ corresponding to $F\cap T_i$, oriented positively with respect to $\mu_i$ and the orientation of $W$.  Then the degeneracy locus on $T_i$ can be expressed as 
$$\delta_i(\Phi_0)= c \mu_i + d \lambda_i$$
where $c$ and $d$ are not necessarily coprime. Note that $c$ must be nonzero. The quotient 
\begin{equation}
\label{equ: fdtc}
c_{T_i}(h) = \frac{d}{c} \in \mathbb{Q}
\end{equation}
is called the {\it fractional Dehn twist coefficient} of $h$ along $T_i\cap F\subset \partial F$.  In the case that the boundary of $F$ is connected, we will drop the the subscript $T_i$.  Fractional Dehn twist coefficients were originally defined to measure the amount of twisting around  $F\cap T_i\subset \partial F$ required to isotope $h$ to its Nielsen-Thurston representative $\varphi$ (see \cite[\S 3.2]{HKMI}), and is an important notion in studying the tightness of the contact structure supported by the open book $(F, h)$ \cite[Theorem 1.1]{HKMI}. A pseudo-Anosov monodromy $h$ is called {\em right-veering} (resp. {\em left-veering}) if $c_{T_i}(h) > 0$ (resp. $c_{T_i}(h) < 0$) for all $i = 1, \cdots, m$.

\subsection{Fenley's asymptotic circle}
\label{subsec: fenley's uc}
Here we describe Fenley's asymptotic circle associated to a pseudo-Anosov flow on a closed $3$-manifold (\cite{Fen12}).   

Given a pseudo-Anosov flow $\Phi$ on a closed, connected, orientable $3$-manifold $W$, let $\widetilde \Phi$ be the pull-back of $\Phi$ to the universal cover $\widetilde W$ of $W$. 

\begin{thm}
{\rm (\cite[Proposition 4.2]{FM01})} 
\label{thm: orbit space}
The orbit space $\mathcal{O}$ of $\widetilde \Phi$ is homeomorphic to $\mathbb R^2$. Moreover, the projection $\pi: \widetilde W \to \mathcal{O}$ is a locally-trivial fibre bundle whose flow line fibres are homeomorphic to $\mathbb R$. 
\qed
\end{thm}
An immediate consequence of the theorem is that closed manifolds, as above, which admit pseudo-Anosov flows are irreducible with infinite fundamental groups, and therefore aspherical. 

The action of $\pi_1(W)$ on $\widetilde W$ descends to one on $\mathcal{O}$. Since the flow lines in $\widetilde W$ inherit a coherent $\pi_1(W)$-invariant orientation, the action of $\pi_1(W)$ on $\mathcal{O}$ is by orientation-preserving homeomorphisms, so we obtain a homomorphism
$$\psi: \pi_1(W) \to \mbox{Homeo}_+(\mathcal{O})$$
Fenley has constructed an ideal boundary for $\mathcal{O}$ over which this action extends.

\begin{thm} 
{\rm (\cite[Theorem A]{Fen12})}
\label{thm: fenley's univ circle}
There is a natural compactification $\mathcal{D} = \mathcal{O} \cup \partial \mathcal{O}$ of $\mathcal{O}$ where $\mathcal{D}$ is homeomorphic to a disk with boundary circle $\partial \mathcal{O}$. The action of $\pi_1(W)$ on $\mathcal{O}$ extends to one on $\mathcal{D}$ by homeomorphisms. 
\qed
\end{thm}
It follows from Fenley's construction that the action of $\pi_1(W)$ on the ideal boundary $\partial \mathcal{O}$ of $\mathcal{O}$ is faithful. That is, the associated homomorphism
$$\rho_\Phi: \pi_1(W) \to \mbox{Homeo}_+(\partial \mathcal{O})$$ 
is injective. We think of $\rho_\Phi$ as taking values in $\mbox{Homeo}_+(S^1)$.

\subsection{Rotation and translation numbers}
\label{subsec: rot and trans numbers}

In this section, we briefly review the translation numbers of elements in $\tG$ and rotation numbers of elements in  $\G$. For details, see \cite[\S 5]{Ghys01}. Though these notions are used implicitly in the proof of Theorem \ref{thm: pa implies co orbifold 2} in \S \ref{subsec: pa flow action on the tree}, the understanding of translation and rotation numbers only becomes essential in \S \ref{subsec: app slope detection}, where we discuss the applications of Theorem \ref{thm: pa implies co orbifold 2} to slope detection.   

Recall that $\tG$ denotes the group of homeomorphisms of the real line which commute with translation by $1$ and can be identified with the universal covering group of $\G$. Given an element $h\in \tG$, the {\em translation number} of $h$, denoted by $\tau(h)$, is defined to be the limit 
\begin{displaymath}
 \lim_{n\to \infty} \frac{h^n(0)}{n}
\end{displaymath}

In particular, if $h$ is translation by $r \in \mathbb R$, i.e., $h(x) = x + r$ for $x\in \mathbb{R}$, then it is easy to verify that $\tau(h) = r$. The following lemma lists some basic properties of the translation number that we will use. See \cite[\S 5]{Ghys01} for proofs. 

\begin{lemma}
\label{lemma: rot and trans}
Let $\tau: \tG \rightarrow \mathbb{R}$ denote the translation number. 
\vspace{-.2cm} 
\begin{enumerate}[leftmargin=*] 
\setlength\itemsep{0.3em}
    \item[{\rm (1)}]  $\tau$ is a homomorphism when restricted to a $\mathbb{Z}\oplus \mathbb{Z}$ subgroup of $\tG$. In particular, $\tau(h^n) = n \tau(h)$ for any $n\in \mathbb{Z}$.
  \item[{\rm (2)}]   $h\in \tG$ has a fixed point on $\mathbb{R}$ if and only if $\tau(h) = 0$.
    \item[{\rm (3)}]   The translation number is invariant under conjugation in $\tG$.
\end{enumerate}

\end{lemma}
Let $f$ be an element in  $\G$ and $\tilde f$ a lift of $f$ in $\tG$. We define the {\it rotation number} of $f$ to be the image of $\tau(\tilde{f})$ in $\mathbb{R}/\mathbb{Z}$. Since any two lifts of $f$ in $\tG$ differ by a translation by an integer, the rotation number of $f$ is well-defined. If $f$ is rotation by $\frac{2\pi a}{n}$, then the rotation number of $f$ is $\frac{a}{n}$ (mod $\mathbb{Z}$).

\section{Asymptotic circle representations of flows on \texorpdfstring{$3$}{3}-orbifolds}
\label{sec: univ circs and locbcs hyp}
Suppose that $\Phi_0$ is a pseudo-Anosov flow on the complement of a link $L = K_1 \cup K_2 \cup \cdots \cup K_m$ in a closed, connected, orientable $3$-manifold $W$ and $\delta_i(\Phi_0)$ its degeneracy locus on $T_i$ for $i = 1, \cdots, m$. Given slopes $\alpha_i$ on $T_i$, it is known that if $|\delta_i(\Phi_0) \cdot \alpha_i| \geq 2$ for each $i$, then $\Phi_0$ extends to a pseudo-Anosov flow on the Dehn filled manifold $X(L)(\alpha_1, \cdots, \alpha_m)$ in such a way that the core of each filling solid torus is a periodic orbit \cite{Fried83}. Then by our discussion in \S \ref{subsec: fenley's uc}, there is a $\G$-representation of $\pi_1(X(L)(\alpha_1, \cdots, \alpha_m))$. In this section, we extend this result to orbifold fillings of $X(L)$ and obtain a $\G$-representation of the orbifold group (cf. Corollary \ref{cor: branched covers of pseudo-Anosov links}). See \cite[Chapter 13]{Thurston77} and \cite[Chapter 2]{BMP03} for discussions of orbifolds and their fundamental groups. The basics of orbifold Dehn filling are covered in \cite[Section 5 of Chapter 2]{BMP03}. 

The existence results of Proposition \ref{prop: pa implies co orbifold} and Corollary \ref{cor: branched covers of pseudo-Anosov links} are used in the proof of our main theorem, Theorem \ref{thm: pa implies co orbifold 2}, to deal with a special case when the leaf space of a certain lamination is degenerate.

\subsection{Well-adapted flows on orbifolds}
\label{subsec: well adapted flows}
We consider $\mathcal{M}$ a closed, connected, oriented $3$-orbifold and $B = B_1 \cup \cdots \cup B_m$ an oriented link in the underlying $3$-manifold $|\mathcal{M}|$ for which the singular set of $\mathcal{M}$ is a union of components of $B$. For each $i$, let $\beta_i$ be the positively oriented meridian of $N(B_i)$, which we also use to denote the associated class in $\pi_1(\mathcal{M})$, and $n_i \geq 1$ be the order of the isotropy along $B_i$.  

Alternatively, one can view such a $3$-orbifold as the result of orbifold filling on a link exterior $X(L)$ along slopes $\alpha_i$. From this viewpoint, $B_i$ is the core of the filling solid torus attached to $T_i$. We denote the resulting orbifold by $X(L)(\alpha_*; n_*)$, where $n_i$ is the order of the isotropy group over $B_i$, $n_* = (n_1, \cdots, n_m)$, and $\beta_i = \alpha_i$ as a slope on $\partial X(L)$. 

We say that a flow $\Phi$ on the underlying $3$-manifold $|\mathcal{M}|$ is {\it well-adapted} to the pair $(\mathcal{M}, B)$ if the following three conditions are satisfied:
\begin{enumerate}
\setlength\itemsep{0.3em}
\item each $B_i$ is an orbit of $\Phi$ and the orientation on $B_i$ agrees with that of the flow; 

\item the restriction $\Phi_0$ of $\Phi$ to the complement $C(B) = |\mathcal{M}|\setminus B$ is pseudo-Anosov; 

\item $n_i|\delta_i(\Phi_0) \cdot \beta_i| \geq 2$ for each $i$, where $\delta_i(\Phi_0)$ is the degeneracy locus of $\Phi_0$ on $\partial N(B     _i)$. 
\end{enumerate}
Let $\Phi_0$ be a pseudo-Anosov flow on the complement of a link $L$ in an oriented,  closed,  connected $3$-manifold $W$, denoted by $C(L)$ as before. We identify $C(L)$ with the interior of $X(L)$. Then there is a flow $\Phi^*$ on $X(L)$ whose restriction to the interior of $X(L)$ is $\Phi_0$. In particular, on a regular neighborhood of any component of $\partial X(L)$, the flow $\Phi^*$ is modeled by the suspension flow of a pseudo-Anosov homeomorphism of a surface with nonempty boundary. 

For each $i$, let $\alpha_i$ be a slope on $T_i$ and $n_i$ a positive integer satisfying $n_i|\delta_i(\Phi_0) \cdot \alpha_i| \geq 2$. Since $|\delta_i(\Phi_0) \cdot \alpha_i|\neq 0$, the linear foliation of $T_i$ by simple closed curves of slope $\alpha_i$ can be isotoped to be everywhere transverse to $\Phi^*|_{T_i}$. Then following \cite[\S 1]{Fried83}, we define a quotient map $\pi: X(L) \rightarrow X(L)(\alpha_*)$ by collapsing every leaf of the linear foliation on $\partial X(L) = \cup_i T_i$ to a point. Then $\Phi^*$ induces a flow $\Phi$ on the Dehn filled manifold $X(L)(\alpha_*)$ and the cores of the filling solid tori are closed orbits of the flow. If $|\delta_i(\Phi_0) \cdot \alpha_i| \geq 2$ for each $i$, the induced flow $\Phi$ on $X(\alpha_*)$ is pseudo-Anosov \cite{Fried83}, though not necessarily otherwise. 

We consider the orbifold $X(L)(\alpha_*; n_*)$, where $\alpha_*= (\alpha_1, \cdots, \alpha_m)$ and $n=(n_1, \cdots, n_m)$. As noted above, the singular set $B=(B_1, \cdots, B_m)$ of $X(L)(\alpha_*; n_*)$ is the union of the cores of the filling solid tori with orientation inherited from the flow $\Phi$ and the underlying manifold of $X(L)(\alpha_*; n_*)$ is $X(\alpha_*)$. It follows that the flow $\Phi$ constructed above is well-adapted to the pair $(X(L)(\alpha_*; n_*), B)$. 

\subsection{\texorpdfstring{$\G$}{Homeo+(S1)}-representations of orbifold groups}
\label{subsec: rep of orbifold groups}
A closed $3$-manifold that admits a pseudo-Anosov flow must be irreducible. The following lemma proves that open manifolds which admit pseudo-Anosov flows satisfy analogous properties; this is needed in the proof of Proposition \ref{prop: pa implies co orbifold} to guarantee that the orbifold $X(L)(\alpha_*; n_*)$ is finitely covered by a manifold. We defer the proof of the lemma to the end of this subsection.

\begin{lemma}
\label{lemma: basic top consequences of cpaf}
Suppose that $\Phi_0$ is a pseudo-Anosov flow on the complement of a link $L$ in a closed, connected, orientable $3$-manifold $W$. Then, 
\begin{enumerate}[leftmargin=*]
\setlength\itemsep{0.3em}
    \item[{\rm (1)}] The exterior $X(L)$ of $L$ in $W$ is irreducible, boundary-incompressible and aspherical.
    \item[{\rm (2)}] If $A$ is an essential annulus in $X(L)$ and $\alpha$ is the slope of a boundary component of $A$ on a torus $T \subset \partial X(L)$, then $|\delta_T(\Phi_0) \cdot \alpha| = 0$.
\end{enumerate}
 
\end{lemma}
\begin{remark}
\label{rem: primeness of links}
Since the exterior of a composite link contains an essential annulus whose boundary components are meridians of some component of the link, part (2) of the lemma implies that if the complement of a link $L$ admits a pseudo-Anosov flow with no meridional degeneracy loci then $L$ is prime.
\end{remark} 

Here is a consequence of Lemma \ref{lemma: basic top consequences of cpaf} which is crucial to our construction of representations of orbifold groups. Before stating it, recall that a {\it teardrop} is a $2$-orbifold of the form $S^2(a)$, where $a > 1$. A {\it spindle} is a $2$-orbifold of the form $S^2(a, b)$, where $1 < a < b$. A fundamental result states that a compact, connected, orientable $3$-orbifold which contains neither teardrops nor spindles is finitely covered by a manifold \cite[Corollary 1.3]{BLP05}. 

\begin{lemma}
\label{lemma: no bad $2$-orbifolds}
Suppose that $\Phi_0$ is a pseudo-Anosov flow on the complement of a link $L$ in a closed, connected, orientable $3$-manifold $W$ and fix essential simple closed curves $\alpha_i$ on $T_i$ and integers $n_i \geq 1$ satisfying $n_i|\delta_i(\Phi_0) \cdot \alpha_i| \geq 1$. Then the orbifold Dehn filling $X(L)(\alpha_*; n_*)$ of the exterior of $L$ is finitely covered by a manifold. 
\end{lemma}

\begin{proof}
Lemma \ref{lemma: basic top consequences of cpaf} shows that $X(L)(\alpha_*; n_*)$ contains neither teardrops nor spindles. The conclusion then follows from \cite[Corollary 1.3]{BLP05}. 
\end{proof}

\begin{prop}
\label{prop: pa implies co orbifold}
Suppose that $\mathcal{M}$, $B$, $\beta_i$, and $n_i$ are as above and that $\Phi$ is a flow on $|\mathcal{M}|$ that is well-adapted to the pair $(\mathcal{M}, B)$. Then there is a faithful representation $\rho_\Phi$ of the orbifold fundamental group $\pi_1(\mathcal{M})$ of $\mathcal{M}:$
$$\rho_\Phi: \pi_1(\mathcal{M}) \to \mbox{{\rm Homeo}}_+(S^1)$$
Further, $\rho_\Phi(\beta_i)$ is conjugate to rotation by $2 \pi/n_i$ for each $i = 1, 2, \ldots, m$. 
\end{prop}

\begin{proof}
By assumption, the restriction $\Phi_0$ of $\Phi$ to $C(B)$ is pseudo-Anosov, so Lemma \ref{lemma: no bad $2$-orbifolds} implies that $\mathcal{M}$ is finitely covered by a manifold. Consider a commutative diagram of covering maps 
\begin{center} 
\begin{tikzpicture}[scale=0.6]
\node at (8, 6) {$\widetilde{\mathcal{M}}$};
\node at (11.4, 3.5) {$W$};
\node at (8, 1) {$\mathcal{M}$}; 
\node at (7.65, 3.5) {$p$};  
\node at (10, 5.3) {$p'$};
\node at (10, 1.8) {$p_1$}; 
\draw [ ->] (8, 5.4) -- (8,1.5);
\draw [ ->] (8.5, 5.6) --(11,3.9); 
\draw [ ->] (10.9, 3.1) -- (8.6,1.4);
\end{tikzpicture}
\end{center} 
where $p$ and $p'$ are universal covers and $p_1$ is a finite degree regular cover from a manifold $W$ to $\mathcal{M}$. Let $\Phi'$ be the lift of the flow $\Phi$ to $W$. Set $(C' , X', B') = p_1^{-1}(C(B), X(B), B)$ and let $\Phi_0'$ be the restriction of $\Phi'$ to $C'$. 

Given a boundary component $T$ of $X'$, suppose that $p_1(T) = \partial N(B_i)$ and let $\upsilon \in H_1(T)$ be the class of the meridional slope of the component of $(p_1)^{-1}(N(B_i))$ containing $T$. Then $\upsilon$ is mapped to $n_i \beta_i \in H_1(\partial N(B_i))$ by the restriction of $p_1$ to $T$. The degeneracy locus $\delta_T(\Phi'_0)$ of $\Phi'_0$ on $T$ is the inverse image of $\delta_i(\Phi_0)$ and so 
$$|\delta_T(\Phi'_0) \cdot \upsilon| = n_i |\delta_i(\Phi_0) \cdot \beta_i| \geq 2$$
It follows that $\Phi'$ is a pseudo-Anosov flow on $W$ and therefore its lift $\widetilde \Phi$ to $\widetilde{\mathcal{M}}$ has orbit space $\mathcal{O} \cong \mathbb R^2$ (Theorem \ref{thm: orbit space}). Since $\widetilde \Phi$ is also the lift of $\Phi$ to $\widetilde{\mathcal{M}}$, it is invariant under the action of $\pi_1(\mathcal{M})$. Thus there is an induced action of $\pi_1(\mathcal{M})$ on $\mathcal{O}$ and as in the proof of Fenley's theorem (cf. Theorem \ref{thm: fenley's univ circle}), this action extends to an action on $\mathcal{D}$, and therefore on $\partial \mathcal{O}$. Let
$$\rho_\Phi: \pi_1(\mathcal{M}) \to \mbox{Homeo}_+(\partial \mathcal{O})$$
be the associated representation. 

The meridional class $\beta_i$ acts on $\widetilde{\mathcal{M}}$ as a rotation by $2 \pi/n_i$ about a component $\widetilde B_i$ of the inverse image of $B_i$ in $\widetilde{\mathcal{M}}$ in the sense determined by the induced orientations on $\widetilde B_i$ and $\widetilde{\mathcal{M}}$. Hence as $\widetilde B_i$ is a flow line of $\widetilde \Phi$, $\rho_\Phi(\beta_i)$ is conjugate to a rotation by $2\pi/n_i$ with respect to the orientation on $S^1 = \partial \mathcal{O}$ determined by the orientation on the flow lines of $\widetilde \Phi$ and that on $\widetilde{\mathcal{M}}$. In particular, $\rho_\Phi$ is injective on each finite subgroup of $\pi_1(\mathcal{M})$.  

The restriction of $\rho_\Phi$ to $\pi_1(W)$ is Fenley's asymptotic circle representation associated to $\Phi'$, so is injective. It follows that each element of $\mbox{kernel}(\rho_\Phi)$ is of finite order and therefore contained in an isotropy subgroup of $\pi_1(\mathcal{M})$, on which we have just seen that $\rho_\Phi$ is injective. Thus the kernel of $\rho_\Phi$ is trivial. 
\end{proof}

\begin{cor}  
\label{cor: branched covers of pseudo-Anosov links} 
Let $L = K_1 \cup \cdots \cup K_m$ be a link in a closed, connected, oriented $3$-manifold $W$ whose complement admits a pseudo-Anosov flow $\Phi_0$. Fix essential simple closed curves $\alpha_i$ on $T_i$ and integers $n_i \geq 1$ satisfying $n_i|\delta_i(\Phi_0) \cdot \alpha_i| \geq 2$. Then
\vspace{-.2cm} 
\begin{enumerate}[leftmargin=*]
\setlength\itemsep{0.3em}
    \item[{\rm (1)}] There is a faithful representation $\rho: \pi_1(X(L)(\alpha_*; n_*)) \to \mbox{{\rm Homeo}}_+(S^1)$. 
    \item[{\rm (2)}] If $\alpha_i$ is oriented positively with respect to the orientation on $X(L)(\alpha_*)$ and that of the core of the  $\alpha_i$-filling solid torus induced by the flow, then  $\rho(\alpha_i)$ is conjugate to rotation by $2 \pi/n_i$ for each $i = 1, 2, \ldots, m$. 
\end{enumerate}
\end{cor}

\begin{proof}
Set $\mathcal{M} = X(L)(\alpha_*; n_*)$. Let $B = B_1 \cup \cdots \cup B_m$ be the link in $|\mathcal{M}| = X(L)(\alpha_*)$ corresponding to the cores of the $\alpha_i$-surgery solid tori. Since $|\delta_i(\Phi_0) \cdot \alpha_i|  \ne 0$ for each $i$, $\Phi_0$ extends to a flow $\Phi$ on $|\mathcal{M}|$ which is well-adapted to the pair $(\mathcal{M}, B)$ once we orient $B$ compatibly with $\Phi$. The corollary then follows immediately from Proposition \ref{prop: pa implies co orbifold}. 
\end{proof}

We finish this subsection by proving Lemma \ref{lemma: basic top consequences of cpaf}. 

\begin{proof}[Proof of Lemma \ref{lemma: basic top consequences of cpaf}] 
Write $L = K_1 \cup \cdots \cup K_m$ and choose slopes $\alpha_i$ on $T_i$ such that $|\delta_i(\Phi_0) \cdot \alpha_i| \geq 2$ for each $i$. Then $\Phi_0$ extends to a pseudo-Anosov flow $\Phi$ on $W' = X(L)(\alpha_1, \ldots, \alpha_m)$ for which the cores of the filling solid tori are flow lines. Theorem 
\ref{thm: orbit space} then implies  that $X(L)$ is covered by a manifold of the form $P \times \mathbb R$, where $P$ is a subsurface of $\mathbb R^2$. Hence $X(L)$  is irreducible and aspherical. If it were boundary-compressible it would be a solid torus and therefore $W' = X(L)(\alpha_1)$ would not be aspherical, contrary to the fact that it admits a pseudo-Anosov flow. This proves (1). 

For (2), suppose that $A$ is an essential annulus in $X(L)$ and $A_0$ is a boundary component of $A$ of slope $\alpha$ on a torus $T \subset \partial X(L)$. Without loss of generality we can assume that $T = \partial N(K_1)$. 

Suppose $|\delta_1(\Phi_0) \cdot \alpha| \geq 1$. Let $\mathcal{O} = X(L)(\alpha;2)$ be the orbifold obtained by an order $2$ orbifold filling to the boundary component $T$ along slope $\alpha$. That is, $\mathcal{O}$ has underlying space $X(L)(\alpha)$ and singular set the core $B_1$ of the $\alpha$-filling solid torus with isotropy $\mathbb Z/2$. By construction $\mathcal{O}$ contains no spindles, since it has only $\mathbb Z/2$ isotropy. Nor does it have any teardrops, since $X(L)$ is boundary-incompressible. Thus there is a finite regular cover $p: M \to \mathcal {O}$, where $M$ is a manifold by Lemma \ref{lemma: no bad $2$-orbifolds}. 

Note that $C(L)$ can be naturally identified with $\mbox{int} (|\mathcal{O}|) \setminus B_1$. We denote the preimages of $C(L)$ and $X(L)$ under $p$ by $\widetilde{C}(L)$ and $\widetilde{X}(L)$ respectively. Let $\widetilde{T}\subset \partial \widetilde{X}(L)$ be a boundary component that covers $T$ and $\tilde{\alpha}$ the slope on $\widetilde{T}$ that covers $\alpha$. Set $M_0 = \widetilde{X}(L)(\tilde{\alpha}) \subseteq M$.

Let $\widetilde{\Phi}_0$ be the pullback flow of $\Phi_0$ under $p|_{\widetilde{C}(L)} : \widetilde{C}(L) \rightarrow C(L)$. Then since the isotropy group along $B_1$ is $\mathbb{Z}/2$,  it follows that 
$$|\delta_{\widetilde{T}}(\widetilde{\Phi}_0) \cdot \tilde{\alpha}| = 2 |\delta_1(\Phi_0) \cdot \alpha| \geq 2$$ 
Hence $\widetilde{\Phi}_0$ can be extended to a pseudo-Anosov flow on the interior of $M_0$. In particular $M_0$ is irreducible, boundary-incompressible and aspherical by $(1)$.

Let $\widetilde A$ be a component of the inverse image of $A$ and $\widetilde{A}_0$ the boundary component of $\widetilde A$ lying over $A_0$. We can assume that $\widetilde A$ is chosen so that $\widetilde A_0 \subset \widetilde T$. Then $\tilde\alpha$ is the slope of $\widetilde{A}_0$. Since $M_0$ is boundary-incompressible, $\partial \widetilde{A}$ must be contained in $\widetilde{T}$; otherwise, $\widetilde{A}$ is contained in an essential disk properly embedded in $M_0$. Then we also have $\partial A \subset T$. 

Let $\mathcal{S} \cong S^2(2,2)$ be the $2$-orbifold in $\mathcal{O}$ obtained from $A$ by attaching two meridional disks of the $\alpha$-filling solid torus. It follows that $\mathcal{S}$ pulls back to a $2$-sphere $\widetilde S$ in $M$, which is the union of $\widetilde A$ and two meridional disks of the component of the inverse image in $M$ of the $\alpha$-filling torus whose boundary is $\widetilde T$. The cover $p$ restricts to a degree $2$ universal cover $\widetilde S \to \mathcal{S}$. The irreducibility of $M$ implies that $\widetilde S$ bounds a $3$-ball $\widetilde B$ in that manifold. 

Suppose that there is a deck transformation $\gamma$ of the cover $M \to \mathcal{O}$ such that $\gamma(\widetilde B) \cap \widetilde B \ne \emptyset$ but $\gamma(\widetilde B) \ne  \widetilde B$. Then either $\gamma^\varepsilon(\widetilde B) \subset \mbox{int} (\widetilde B)$ for some $\varepsilon \in \{\pm 1\}$ or $M = \widetilde B \cup \gamma(\widetilde B)$. The former is ruled out since $\gamma$ is of finite order, while the latter is ruled out by the asphericity of $M$. Thus $\gamma(\widetilde B) \cap \widetilde B \ne \emptyset$ implies that $\gamma(\widetilde B) = \widetilde B$. It follows that the stabiliser of $\widetilde B$ coincides with that of $\widetilde S$, and is therefore isomorphic to $\mathbb Z/2$. Then we have a $2$-fold cover $\widetilde B \to \mathcal{B} \subset \mathcal{O}$ of orbifolds, which can be thought of as a $2$-fold cover branched over the properly embedded arc $|\mathcal{B}| \cap B_1$ in the underlying space $|\mathcal{B}|$ of $\mathcal{B}$. Since $\widetilde B$ is a $3$-ball, $|\mathcal{B}|$ is a simply-connected compact $3$-manifold with boundary a $2$-sphere, and therefore a $3$-ball. But then as $\widetilde B$ is a $3$-ball, $|\mathcal{B}| \cap B_1$ is unknotted in $|\mathcal{B}|$ by the $\mathbb Z/2$-Smith Conjecture, and therefore $A$ is boundary parallel in $X(L)$, which contradicts our assumption that it is essential. This final contradiction shows that $|\delta_T(\Phi_0) \cdot \alpha| = 0$. 
\end{proof}

\section{Recalibrating \texorpdfstring{$\mathbb R$}{R}-order trees and \texorpdfstring{$\mbox{Homeo}_+(S^1)$}{Homeo+(S1)}-representations} 
\label{sec: psA flows and $LO$ cbcs}

\subsection{\texorpdfstring{$\mathbb R$}{R}-order trees} 
\label{subsec: R-order tree}

An $\mathbb R$-order tree is a set $T$ together with a family $\mathcal{S}(T)$ of totally ordered subsets whose elements are called {\it segments}. Each segment $\sigma$ is order isomorphic to a closed interval in the reals whose least and greatest elements (its endpoints) are assumed to be distinct and are denoted $\sigma^-$ and $\sigma^+$ respectively. The {\it inverse} of a segment $\sigma$ is $\sigma$ endowed with the opposite order and is denoted $-\sigma$. It is assumed that the following conditions hold:

\begin{itemize}
\setlength\itemsep{0.2em}
\item $\sigma\in \mathcal{S}(T)$ if and only if $-\sigma\in \mathcal{S}(T)$;  
\item a closed subinterval of a segment with more than one element is a segment;  
\item any two elements of $T$ can be joined by a sequence $\sigma_1, \cdots, \sigma_k \in \mathcal{S}(T)$ with $\sigma_i^{+} = \sigma_{i+1}^-$ for $i = 1, \cdots, k-1$;  
\item if $\sigma_0\sigma_1\cdots \sigma_{k-1}$ is a cyclic word of segments with $\sigma_i^{+} = \sigma_{i+1}^-$ for all $i$ (mod $k$), then each $\sigma_j$ can be written as a concatenation of subsegments yielding a cyclic word $\rho_0\rho_1\cdots\rho_{n-1}$ which becomes the trivial word when adjacent inverse segments are canceled; 
\item if $\sigma_1$ and $\sigma_2$ are segments whose intersection is $\sigma_1^+ = \sigma_2^-$,  then $\sigma_1\cup \sigma_2\in \mathcal{S}(T)$;
\item $T$ is the union of countably many segments.
\end{itemize}
We will sometimes write $[\sigma^-, \sigma^+]$ for $\sigma$, $(\sigma^-, \sigma^+]$ for $\sigma \setminus \sigma^-$, etc. An $\mathbb R$-order tree is endowed with the weak topology with respect to its segments. 

It is interesting to contrast this definition with that of a real tree \cite[Section 2]{Bst02}, which is a metric space satisfying similar conditions, though its segments are required to be isometric to closed intervals. There is no natural metric on an $\mathbb R$-order tree. 

Order trees were defined in \cite{GO89} to model the structure of the leaf space of an essential lamination $\Lambda$ on $\mathbb R^3$. The definition above was taken from \cite{GK97}. (See also \cite{RS01}.)  

Let $\Lambda$ be an essential lamination on a closed $3$-manifold $M$. We may assume that $\Lambda$ has no isolated leaves. Otherwise, one can perform the well-known ``blowup'' operation on isolated leaves \cite[Operation 2.1.1]{Gab92}. If $\Lambda$ is the stable or unstable lamination of a pseudo-Anosov flow, then there are no isolated leaves. Let $\widetilde{\Lambda}$ denote the pullback of $\Lambda$ to the universal cover $\widetilde{M} \cong \mathbb{R}^3$. Then the leaf space of $\widetilde{\Lambda}$, denoted by $T(\widetilde{\Lambda})$, is defined to be the quotient $\widetilde{M}/\hspace{-1.5mm} \sim$ where $x\sim y$ if $x$ and $y$ are on the same leaf of $\widetilde{\Lambda}$ or they are in the same complementary region of $\widetilde{\Lambda}$.  

The leaf space $T(\widetilde{\Lambda})$ is an $\mathbb{R}$-order tree \cite{GK97}. The segments of $T(\widetilde{\Lambda})$ are the images of transverse arcs to $\widetilde{\Lambda}$ under the natural map $v: \mathbb{R}^3 \rightarrow T(\widetilde{\Lambda})$. One of the main results of \cite{GK97} states that this association induces a one-to-one correspondence between homeomorphism classes of essential laminations on $\mathbb R^3$ and isomorphism classes of $\mathbb R$-order trees (\cite[Theorem 8.1]{GK97}). 

Analogous definitions and results hold for essential laminations on $\mathbb R^2$, though here the lamination's embedding in the plane endows the tree with extra structure in the form of local circular orderings, which we will discuss in the next section. Theorem 4.2 of \cite{GK97} states that every $\mathbb R$-order tree endowed with this extra structure is isomorphic to the leaf space of some essential lamination on $\mathbb R^2$, while \cite[Theorem 6.9]{GK97} states that two essential laminations on $\mathbb R^2$ with isomorphic $\mathbb R$-order trees are equivalent up to splitting along leaves.

\subsection{Cyclically ordered \texorpdfstring{$\mathbb R$}{R}-order trees}
\label{subsec: cyclically ordered R-order tree}
A {\it circular order} on a set $E$  is a function $c:E^3 \rightarrow \{ -1, 0, 1\}$ satisfying:
\begin{itemize}
    \setlength\itemsep{0.2em}
\item if $e_1, e_2, e_3 \in E$, then $c(e_1, e_2, e_3) = 0$ if and only if $e_i = e_j$ for some $i \ne j$;
\item the {\it cocyle condition} : for all $e_1, e_2, e_3, e_4 \in E$ we have
$c(e_2, e_3, e_4) - c(e_1, e_3, e_4) + c(e_1, e_2, e_4) - c(e_1, e_2, e_3) = 0$.
\end{itemize}
Circular orders abstract the positional arrangement of points on an oriented circle, where the condition $c(e_1, e_2, e_3) = 1$ means that 
$e_1, e_2, e_3$ are distinct points and the oriented segment from $e_1$ to $e_3$ passes through $e_2$. Countable circularly ordered sets can be embedded in the circle in an order preserving way. 

\begin{figure}[ht]
    \includegraphics[scale=0.7]{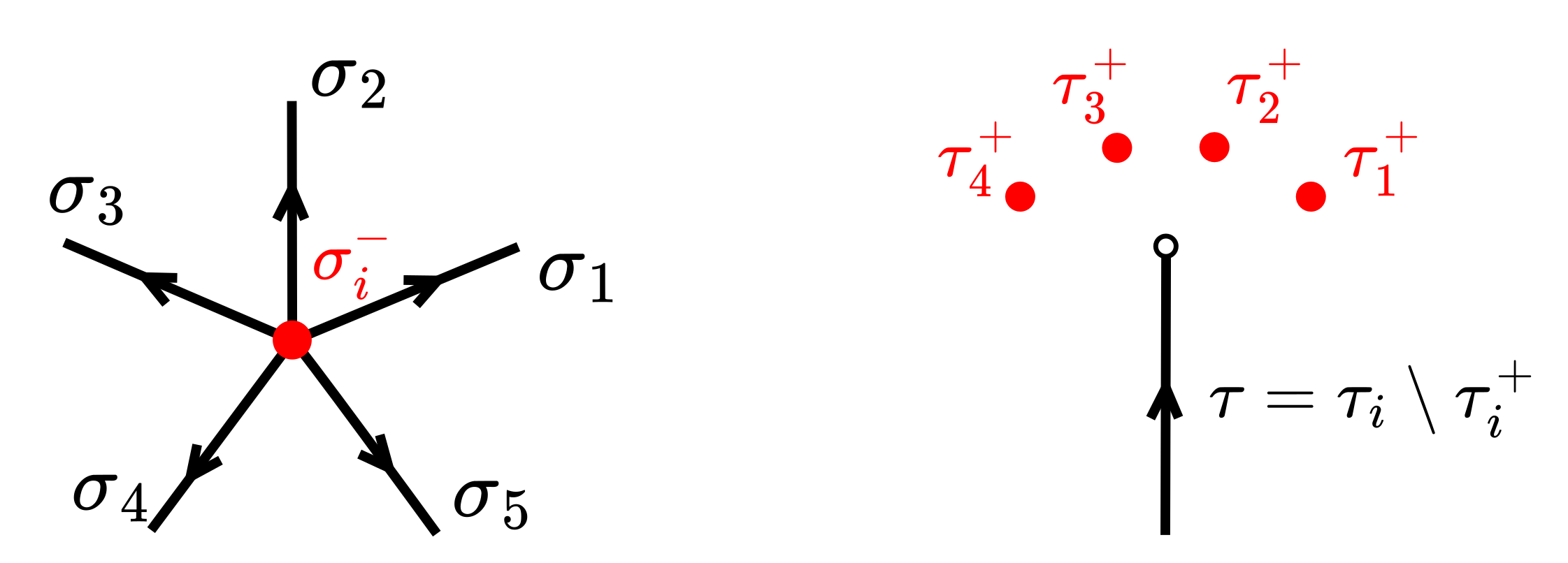}
    \caption{Two types of local circular orders on an $\mathbb{R}$-order tree}
    \label{fig: circular order tree}
\end{figure}

An immediate consequence of these conditions is that if $e_1, e_2, e_3 \in E$ and $\tau \in S_3$ a permutation, then 
\begin{equation} 
\label{eqn: 3-permutation}
c(e_{\tau(1)}, e_{\tau(2)}, e_{\tau(3)})) = \mbox{sign}(\tau)c(e_1, e_2, e_3)  
\end{equation}  
Consider an $\mathbb R$-order tree $T$ with $\mathcal{S}(T)$ its family of ordered segments. A {\it cataclysm} of $T$ is a subset of $T$ with cardinality two or more of the form $\overline{\tau} - (\tau \setminus \tau^+)$, where $\tau$ is a segment of $T$ (cf. \cite[\S 3.4]{CD03}). For instance, in Figure \ref{fig: circular order tree}, the set $\{\tau_1^+, \tau_2^+, \tau_3^+, \tau_4^+\}$ is contained in a cataclysm. Following \cite[\S 2]{GK97}, a {\it cyclic ordering} on $T$ consists of two sets of local circular orders: 
\begin{itemize}
\setlength\itemsep{0.2em}
\item[(a)] a circular order on any finite set of segments $\{\sigma_i\}$ such that $\sigma_i^- =\sigma_j^-$ for all $i, j$ and $(\sigma_i \setminus \sigma_i^-)\cap (\sigma_j \setminus \sigma_j^-) =\emptyset$ for all $i\neq j$. Further, the circular order on a subset of $\{\sigma_i\}$ or on a set of initial subsegments of $\sigma_i$'s is that induced from $\{\sigma_i\}$.  In this case we call the common point $\sigma_i^-$ a {\it vertex} of $T$.
\item[(b)] a circular order on any finite set of segments $\{\tau\} \cup \{\tau_i\}_{i=1}^n$, where $n \geq 2$, $\tau_i^+\neq \tau_j^+$ for $i\neq j$ and $\tau = \tau_i \setminus \tau_i^+ = \tau_j \setminus \tau_j^+$ for all $i, j$. Further, the circular order on a subset of $\{\tau\} \cup \{\tau_i\}$ or on a set of terminal subsegments of $\tau$ and $\tau_i$'s is that induced from $\{\tau\} \cup \{\tau_i\}$. In this case the set $\{\tau_i^+\}_{i=1}^n$ is contained in a cataclysm of $T$.
\end{itemize} 

Thus a cyclic order on $T$ determines circular orders on the germs of segments incident to vertices and cataclysms. 
 
Following Roberts and Stein \cite{RS01}, we define a {\it cusp} of $T$ to be an equivalence class of pairs of segments $(\sigma, \tau)$ with $\sigma \setminus \sigma^+ = \tau \setminus \tau^+$ but $\sigma^+ \neq \tau^+$, where we say that $(\sigma_1, \tau_1)$ is equivalent to $(\sigma_2, \tau_2)$ if $\sigma_1^+ = \sigma_2^+$, $\tau_1^+ = \tau_2^+$ and there is an $x \in (\sigma_1 \setminus \sigma_1^+) \cap (\sigma_2 \setminus \sigma_2^+)$ such that $\sigma_i = [\sigma_i^-, x] \cup [x, \sigma_i^+]$ and $\tau_i = [\tau_i^-, x] \cup [x, \tau_i^+]$ for both values of $i$. A cusp is determined uniquely by the pair $\{\sigma^+, \tau^+\}$, which we use to denote it. The reader will verify that if $\{p, q\}$ is a cusp, then $p$ and $q$ are contained in a cataclysm.  A {\it stem} of a cusp $\{p, q\}$ is any of the intervals $[\sigma^-, \sigma^+)$, where $(\sigma, \tau)$ represents $\{p, q\}$. 

\begin{figure}[ht]
    \centering
    \includegraphics[scale=0.7]{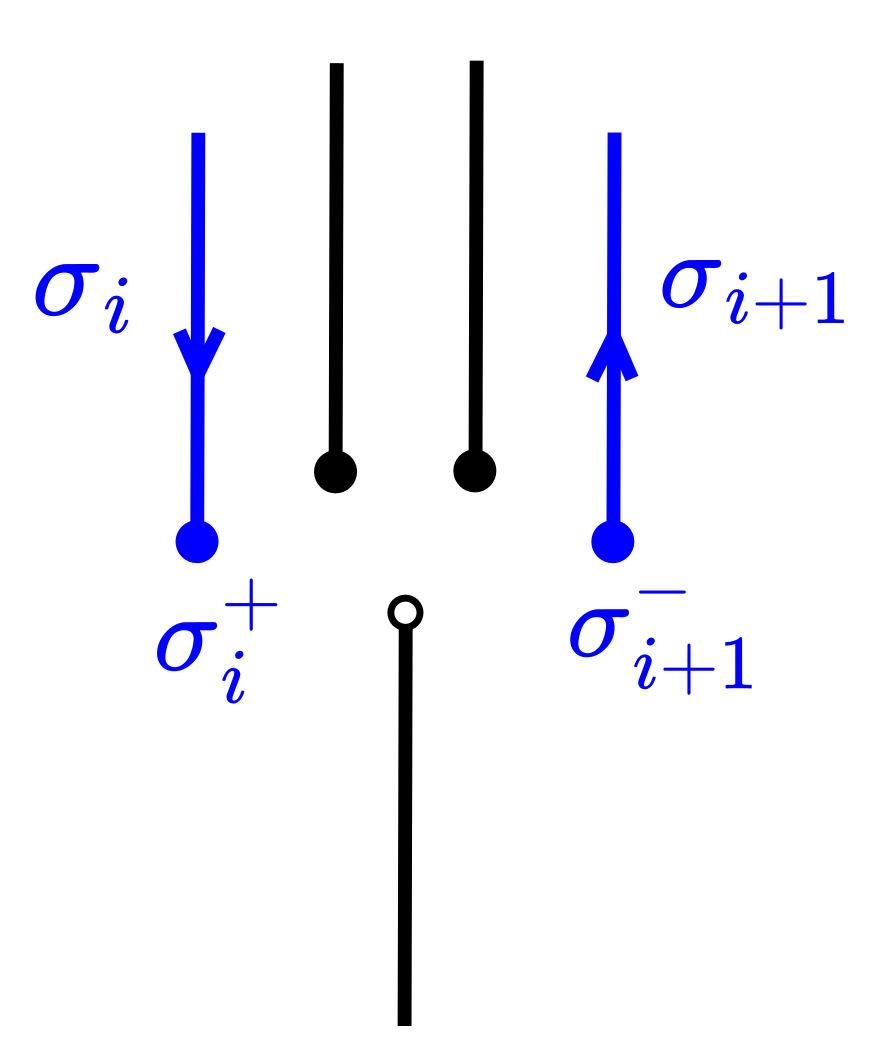}
    \caption{The segments colored in blue can be a part of a geodesic spine}
    \label{fig: geodesic spine}
\end{figure}

The {\it geodesic spine} between points $x, y$ of $T$ (\cite[Definition 3.5, Theorem 3.6]{RS01}) is the intersection of the images of all paths between them. Moreover, it is easy to use the amalgamation of segments to see that the geodesic spine between $x$ and $y$ can be expressed uniquely as a union $\sigma_1 \cup  \sigma_2 \cup  \cdots \cup \sigma_n$ of segments where $\sigma_1^- = x, \sigma_n^+ = y$ and $\{\sigma_i^+, \sigma_{i+1}^-\}$ is a cusp for $1 \leq i \leq n-1$. See Figure \ref{fig: geodesic spine}. It is possible that $\sigma_1$ and/\hspace{-.2mm}or $\sigma_n$ are the degenerate segments $x$ and/\hspace{-.2mm}or $y$ when the latter belong to a cusp. 

More generally, a {\it geodesic spine} is a subset of $T$ which can be expressed as a nested union of geodesic spines between pairs of points.

For each triple of distinct points $(x, y, z)$ of an $\mathbb R$-order tree, Roberts and Stein defined a subset $Y_{(x, y, z)}$ of $T$ representing the unique place in $T$ at which the geodesic spine from $y$ to $x$ diverges from that from $y$ to $z$ (\cite[page 182]{RS01}). It is either a point or a cusp, where the latter occurs if and only if both geodesic spines contain a stem of the cusp and $x$ and $z$ lie in different components of the complement  of one (and hence all) of the cusp's stems in $T$. 

A {\it ray} based at $x \in T$ is a proper embedding $r: ([0, \infty), 0) \to (T, x)$. By definition, a ray can only go through at most one point in each cataclysm; otherwise, the map $r$ cannot be an embedding. One can see this in Figure \ref{fig: geodesic spine}, where any continuous path going through both blue segments must travel downward into the stem of the cusp and then backtrack upward.  

We consider two rays $r_1, r_2$ to be equivalent, written $r_1 \sim r_2$, if there are real numbers $t_1, t_2 \geq 0$ such that $r_1([t_1, \infty)) = r_2([t_2, \infty))$, and define the set of ends of $T$ to be  
$$\mathcal{E}(T) = \{\mbox{rays in } T \}/\hspace{-1mm}\sim$$
Given $e \in \mathcal{E}(T)$ and point $x \in T$, there is a (unique) geodesic spine in $T$ between $x$ and $e$ obtained by concatenating a geodesic spine between $x$ and some point $y$ of $T$ with a ray based at $y$ in the class of $e$.  

The following is stated in the proof \cite[Theorem 3.8]{CD03}. For completeness we include a proof.  

\begin{lemma}
\label{lem: end co}
The set of ends $\mathcal{E}(T)$ of a cyclically ordered $\mathbb R$-order tree $T$ admits a natural circular ordering.
\end{lemma}

\begin{proof}
To define a circular ordering $c$ on $\mathcal{E}(T)$ associated to the cyclic ordering on $T$, let $e_* = (e_1, e_2, e_3) \in \mathcal{E}(T)^3$ and set $c(e_*) = 0$ if $e_i = e_j$ for some $i \ne j$. 

If the $e_i$ are distinct, choose representative rays $r_i$ for the $e_i$ and increasing sequences $\{x_n\}$ in $r_1$, $\{y_n\}$ in $r_2$, and $\{z_n\}$ in $r_3$ which limit, respectively, to $e_1, e_2, e_3$. Our hypotheses imply that for large $n$ the points $x_n, y_n, z_n$ are distinct and the subsets $Y(x_n, y_n, z_n)$, $Y(y_n, z_n, x_n)$, and $Y(z_n, x_n, y_n)$ stabilise. By taking subsequences we can assume that $Y(x_n, y_n, z_n), Y(y_n, z_n, x_n)$, and $Y(z_n, x_n, y_n)$ are independent of $n$. 

Theorem 3.10 of \cite{RS01} implies that there are segments $\sigma_1, \sigma_2, \sigma_3$ of $T$ such that $(\sigma_i - \sigma_i^-)\cap (\sigma_j -\sigma_j^-) =\emptyset$ for $i \ne j$ and one of the following situations arises: 

\begin{itemize}
\setlength\itemsep{0.2em}
\item $\sigma_1^- = \sigma_2^- = \sigma_3^- = Y(x_n, y_n, z_n) = Y(z_n, x_n, y_n) = Y(y_n, z_n, x_n)$ is a vertex $v(e_*)$ of $T$. Moreover, $\sigma_i \setminus \sigma_i^-$ is contained in the component of $T \setminus v(e_*)$ containing each $x_n$ when $i = 1$, $y_n$ when $i = 2$, and $z_n$ when $i = 3$. See Figure \ref{fig: Y case 1}, where the components of $T\setminus v(e_*)$ containing $\{x_n\}$, $\{y_n\}$, $\{z_n\}$ are illustrated. (Note that these points are not necessarily on $\sigma_i$.)  

\begin{figure}[ht]
    \includegraphics[scale=0.7]{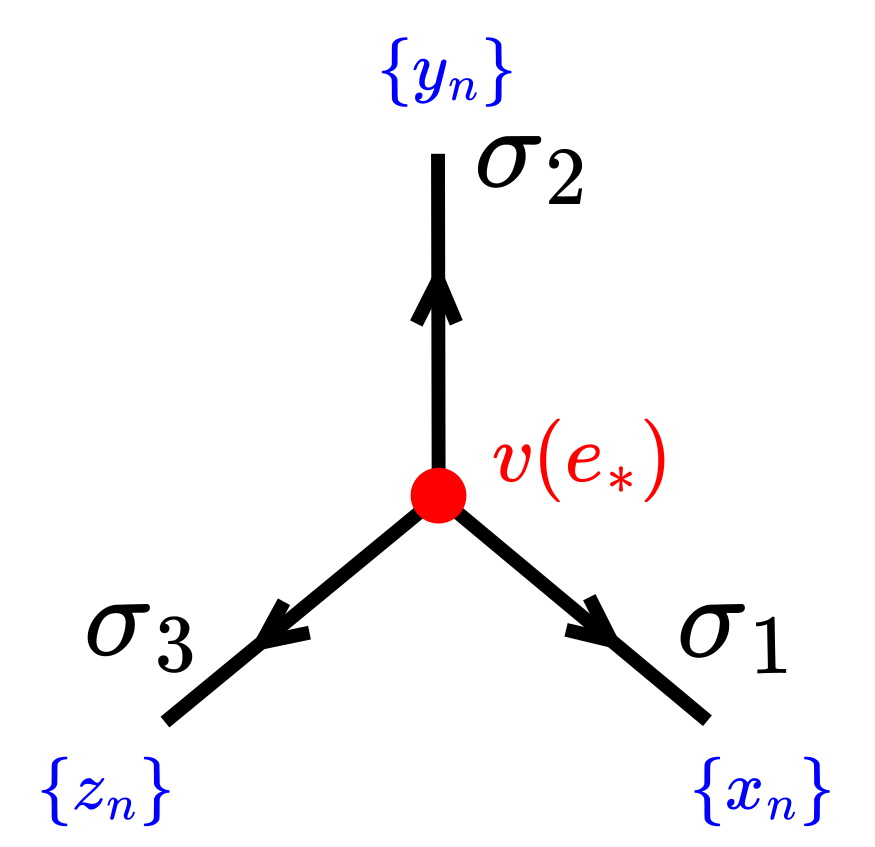}
    \caption{}
    \label{fig: Y case 1}
\end{figure}

 \item $\sigma_1^- = Y(z_n, x_n, y_n), \; \sigma_2^- = Y(x_n, y_n, z_n), \; \sigma_3^- = Y(y_n, z_n, x_n)$ are distinct points contained in a cataclysm of $T$, also denoted $v(e_*)$. Moreover, $\sigma_i \setminus \sigma_i^-$ is contained in the component of $T \setminus v(e_*)$ containing each $x_n$ when $i = 1$, $y_n$ when $i = 2$, and $z_n$ when $i = 3$. (Figure \ref{fig: Y case 2})
 
 \begin{figure}[ht]
    \includegraphics[scale=0.7]{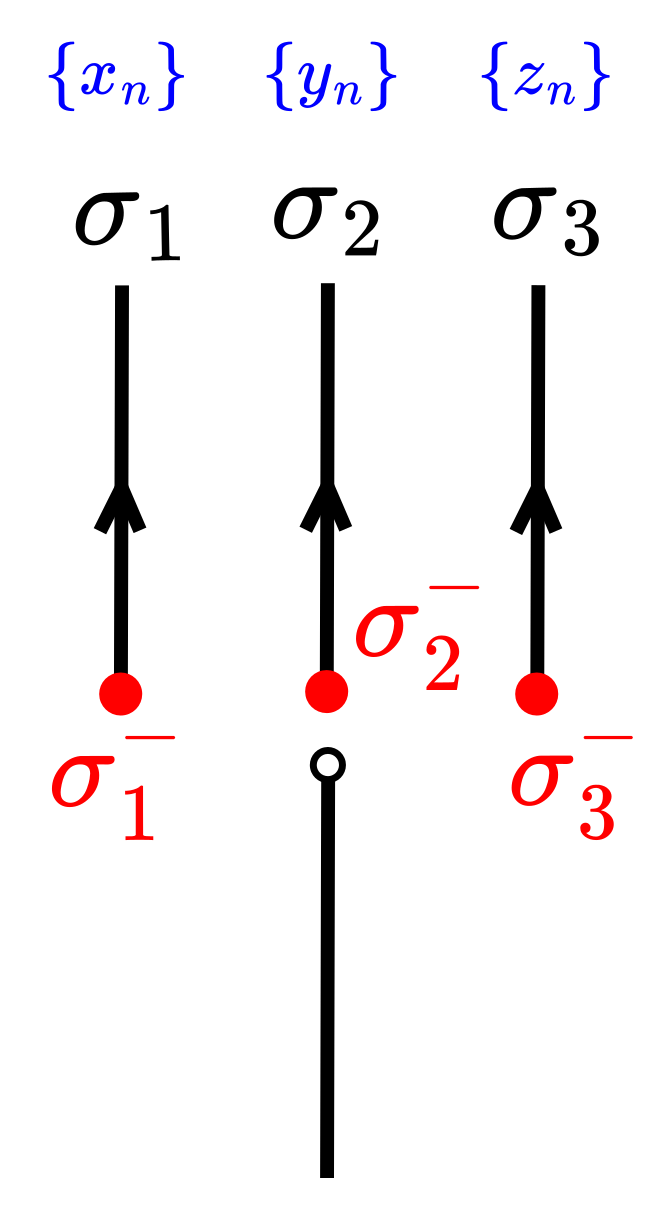}
    \caption{}
    \label{fig: Y case 2}
\end{figure}

\item up to switching $e_3$ and one of $e_1, e_2$: $\sigma_1^- = Y(z_n, x_n, y_n)$ and $\sigma_2^- = Y(x_n, y_n, z_n)$ are distinct points contained in a cataclysm $v(e_*)$ of $T$, and $\{\sigma_1^-, \sigma_2^-\} = Y(y_n, z_n, x_n)$. Moreover, $\sigma_i \setminus \sigma_i^-$ is contained in the component of $T \setminus v(e_*)$ containing each $x_n$ when $i = 1$, $y_n$ when $i = 2$, and $z_n$ when $i = 3$. (Figure \ref{fig: Y case 3})

\begin{figure}[ht]
    \includegraphics[scale=0.7]{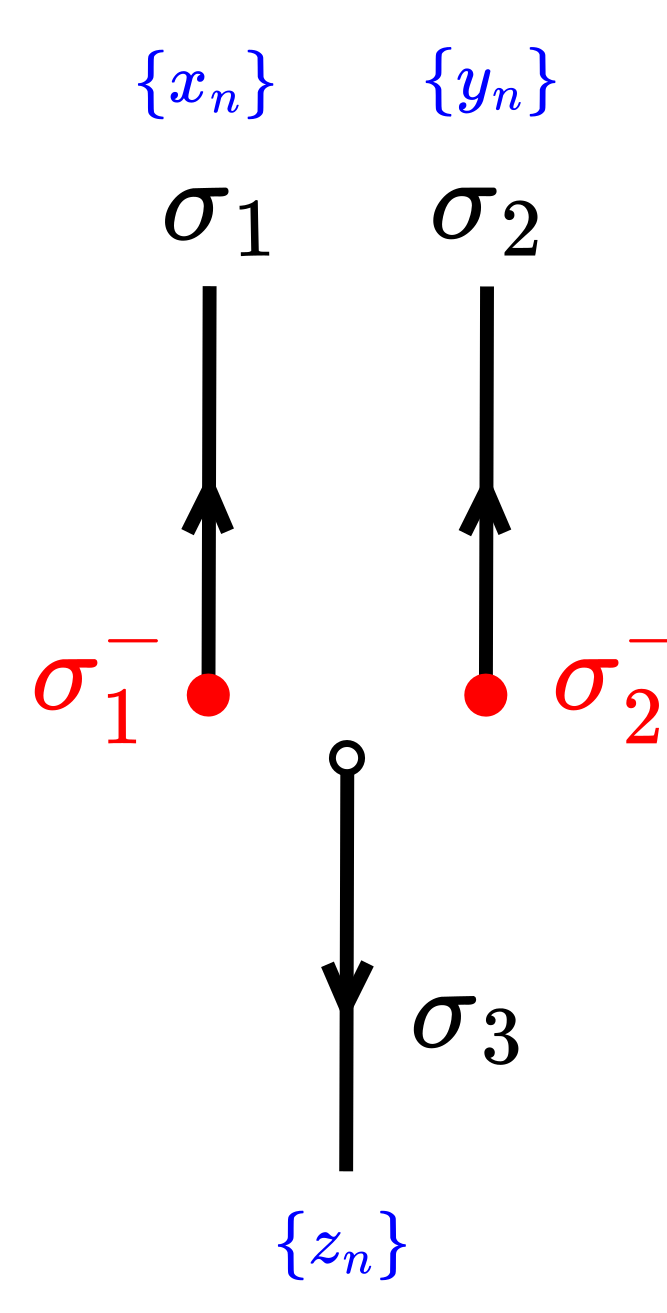}
    \caption{}
    \label{fig: Y case 3}
\end{figure}
\end{itemize}

In each case the cyclic ordering on $T$ determines a circular ordering $c_{v(e_*)}$ on the germs of segments with an endpoint in $v(e_*)$. Set
\begin{equation}
    \label{equ: circular order end}
    c(e_1, e_2, e_3) = c_{v_{(e_*)}}(\sigma_1, \sigma_2, \sigma_3)
\end{equation}
To complete the proof we must show that 
\begin{equation}
    \label{equ: circular order end co cycle}
    s = c(e_2, e_3, e_4) - c(e_1, e_3, e_4) + c(e_1, e_2, e_4) - c(e_1, e_2, e_3)
\end{equation}
is zero for each choice of $e_1, e_2, e_3, e_4 \in \mathcal{E}(T)$. 

This is readily verified if $e_i = e_j$ for some $i \ne j$, since $v(e_*)$ is invariant under permutations of $e_*$. 

Suppose then that $e_1, e_2, e_3, e_4$ are distinct elements of $\mathcal{E}(T)$ and let $e_* = (e_1, e_2, e_3)$.  For each $i = 1, \cdots, 4$, there is a geodesic spine $\gamma_i$ from $v(e_*)$ to $e_i$ obtained by concatenating a geodesic spine $c_i$ from $v(e_*)$ to a point $p_i$ of $T$ with a ray $r_i$ in the class of $e_i$ based at $p_i$. We assume that $c_i$ contains only one point of $v(e_*)$, which prohibits $c_i$ from starting off by jumping between different points in $v(e_*)$ when $v(e_*)$ is a cataclysm. For instance, in Figure \ref{fig: Y case 3}, $\sigma_2^- \cup \sigma_1$ is a geodesic spine that contains two points of $v(e_*)$. 

Let $\sigma_1, \sigma_2, \sigma_3, \sigma_4$ denote initial segments of $\gamma_1, \gamma_2, \gamma_3, \gamma_4$, oriented so that $\sigma_i^- \in v(e_*)$ for each $i$. Then by our discussion of the three cases above (cf. Figures \ref{fig: Y case 1}, \ref{fig: Y case 2} and \ref{fig: Y case 3}), $\sigma_i \setminus \sigma_i^-$ are pairwise disjoint for $i = 1, 2, 3$. If $\sigma_4\setminus\sigma_4^-$ is also disjoint from $\sigma_i \setminus \sigma_i^-$ for $i = 1, 2, 3$, then the identity $s = 0$ follows from the fact that $c_{v(e_*)}$ is a circular ordering. Otherwise, there is a unique $i \in \{1,2,3\}$ such that $(\sigma_4\setminus \sigma_4^-) \cap (\sigma_i\setminus \sigma_i^-) \neq \emptyset$.  
 
Sketching some possible configurations for the $\{\sigma_i\}$ shows intuitively why the cocycle condition (\ref{equ: circular order end co cycle}) holds in this case, though the formal verification is somewhat tedious. As such we depict the case that $i=1$ and $v(e_*)$ is a vertex in Figure \ref{fig: config sigma_i} to help orient the reader through the formal argument. 

\begin{figure}[ht]
    \centering 
    \includegraphics[scale=0.6]{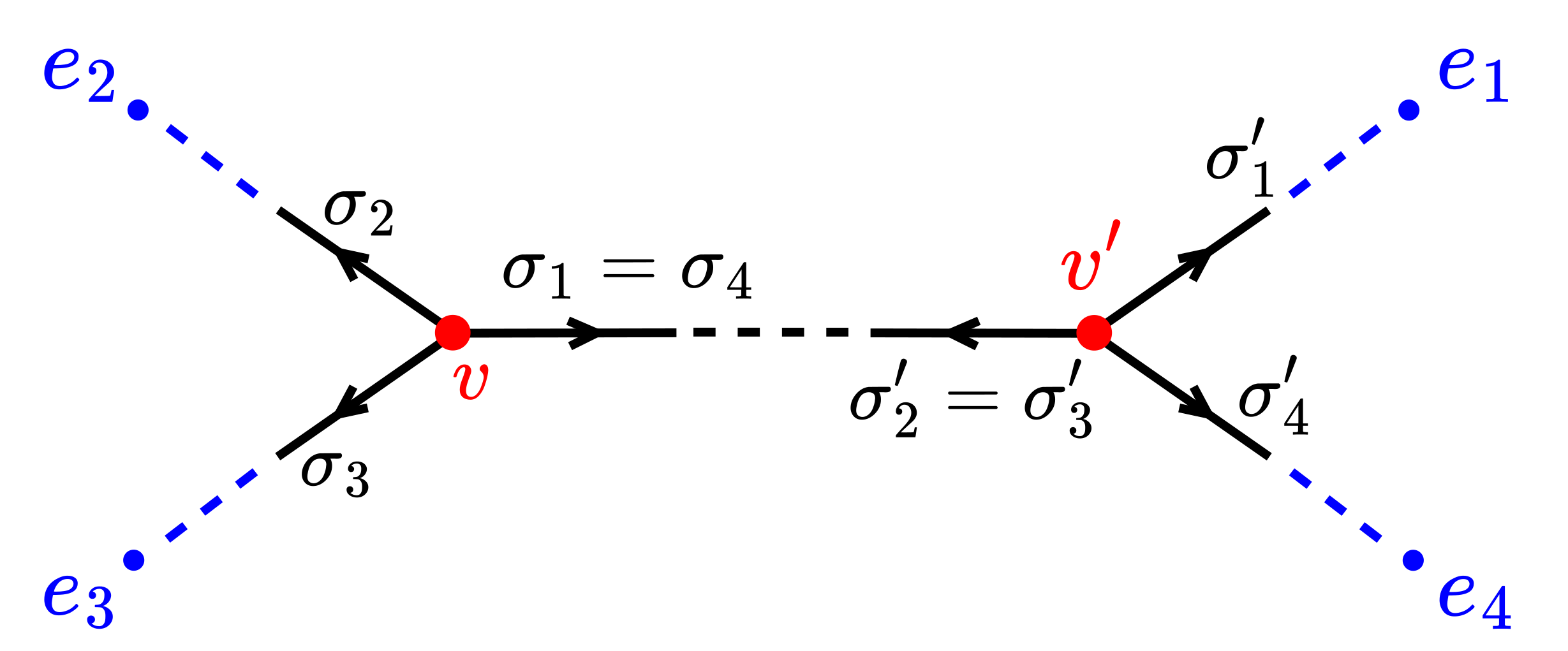}
    \caption{Here, $i = 1$ and we use $v$ to denote $v(e_*)$, $v'$ to denote $v(e_i, e_k, e_4) = v(e_i, e_j, e_4)$ and suppose that $v$ and $v'$ are vertices. In the case that one of them is a cataclysm, the figure is similar (cf. Figure \ref{fig: Y case 2} and Figure \ref{fig: Y case 3}).}
    \label{fig: config sigma_i}
\end{figure}

For the general case, set $\{j, k\} = \{1, 2, 3\} \setminus \{i\}$, in which case the reader will verify that $v(e_i, e_j, e_4) = v(e_i, e_k, e_4)$.  After possibly shrinking $\sigma_4$ and $\sigma_i$, we can suppose that $\sigma_4 = \sigma_i$. Let $\sigma_i', \sigma_4'$, $\sigma_j'$, $\sigma_k'$ be initial segments of the geodesic spines based at $v'$ in the class of $e_i, e_4$, $e_j$ and $e_k$. Note that we can choose $\sigma_j' = \sigma_k'$.  Then, 
\begin{eqnarray} 
s =  \left\{ 
\begin{array}{ll}
c_{v}(\sigma_2, \sigma_3, \sigma_1) - c_{v'}(\sigma_1', \sigma_3', \sigma_4') + c_{v'}(\sigma_1', \sigma_2', \sigma_4')-  c_{v}(\sigma_1, \sigma_2, \sigma_3) & \mbox{ if } i = 1 \\ 
c_{v'}(\sigma_2', \sigma_3', \sigma_4') - c_v(\sigma_1, \sigma_3, \sigma_2) + c_{v'}(\sigma_1', \sigma_2', \sigma_4') - c_v(\sigma_1, \sigma_2, \sigma_3) & \mbox{ if } i = 2 \\
c_{v'}(\sigma_2', \sigma_3', \sigma_4') - c_{v'}(\sigma_1', \sigma_3', \sigma_4') + c_v(\sigma_1, \sigma_2, \sigma_3) - c_v(\sigma_1, \sigma_2, \sigma_3) & \mbox{ if } i = 3 
\end{array} \right. \nonumber 
\end{eqnarray}
It is simple to see that each of these expressions is zero using (\ref{eqn: 3-permutation}) and the fact that we can take $\sigma_2' = \sigma_3'$ for the first, $\sigma_1' = \sigma_3'$ for the second, and $\sigma_1' = \sigma_2'$ for the third.  This completes the proof. 

\end{proof}

\subsection{Group actions on \texorpdfstring{$\mathbb{R}$}{R}-order trees}
We say that a group $G$ acts on an $\mathbb R$-order tree $T$ if it acts on the underlying space preserving $\mathcal{S}(T)$. We say that it acts on a cyclically ordered $\mathbb R$-order tree if it acts on the $\mathbb R$-order tree preserving the local circular orders. Since any such action satisfies
$$Y_{(g \cdot x,g \cdot y, g \cdot z)} = g(Y_{(x,y,z)})$$
for triples of distinct points $x, y, z \in T$, we deduce the following lemma.

\begin{lemma}
\label{lemma: ord pres on ends}
If a group $G$ acts on a cyclically ordered $\mathbb R$-order tree, then $G$ acts on the set of ends of $T$ by order-preserving automorphisms.
\qed
\end{lemma}

A {\it circular order} on a group $G$ is a circular order on the set $G$ which satisfies
$$c(g\cdot g_1, g\cdot g_2, g\cdot g_3) = c(g_1,g_2,g_3)$$
for all $g, g_1, g_2, g_3 \in G$.

An automorphism of a circularly ordered set $(E, c)$ is a bijection of $E$ which preserves $c$. The group of such automorphisms will be denoted by $\mbox{Aut}(E, c)$.  

\begin{lemma}
\label{lemma: aut is co}
If $(E, c)$ is a circularly ordered set, then the group $\mbox{Aut}(E, c)$ is circularly ordered. 
\end{lemma}

\begin{proof}
A proof in the case that $E$ is the circle with its natural circular order is contained in \cite[Theorem 2.2.14]{Cal04}. The general case is similar, though we sketch the argument for use in the proof of Theorem \ref{thm: pa implies co orbifold 2}. 

Set $G = \mbox{Aut}(E, c)$, fix $e \in E$ and let $G_e \leq G$ be the stabiliser of $e$. We obtain a $G$-invariant circular order $c_e$ on the set of cosets $G/G_e$ via the $G$-invariant embedding $gG_e \mapsto g\cdot e$: 
$$c_e(g(g_1 G_e), g(g_2 G_e), g(g_3 G_e)) = c_e(g_1 G_e, g_2 G_e, g_3 G_e) \mbox{ for all } g, g_1, g_2, g_3 \in G$$ 
Hence $G$ is circularly ordered if $G_e = \{1\}$. 

If $G_e \ne \{1\}$, observe that $E \setminus \{e\}$ admits a $G_e$-invariant total order defined by
$$e_1 <e_2 \Leftrightarrow c(e, e_1, e_2) = 1$$
Since $G_e$ acts faithfully on $E \setminus \{e\}$ preserving this total order, it admits a left-order $<_e$ (\cite{Con59}). 
A circular order $c'$ on $G$ can then be obtained from the sequence $1 \to G_e \to G \to G/G_e$ by piecing together the left-order $<_e$ on $G_e$ and the circular order $c_e$ on $G/G_e$: If $g_1, g_2, g_3 \in G$ set 
\vspace{-.2cm} 
\begin{itemize}
\item $c'(g_1, g_2, g_3) =  0$ if $g_i = g_j$ for some $i \ne j$;
\vspace{.2cm} \item $c'(g_1, g_2, g_3) = c_e(g_1G_e, g_2G_e, g_3G_e)$ if $g_1G_e, g_2G_e, g_3G_e$ are distinct cosets. 
\end{itemize}

In the case when $g_1, g_2, g_3$ are distinct, but $g_1G_e, g_2G_e, g_3G_e$ are not, set:

\begin{itemize}
\item $c'(g_1, g_2, g_3) = (-1)^{(j-i) + 1}\mbox{sign}_{<_e}(g_i^{-1} g_j)$  if $g_i G_e = g_j G_e \ne g_k G_e \mbox{ where } i < j$; 
\vspace{.2cm} \item $c'(g_1, g_2, g_3) = \mbox{sign}(\tau)$ if $g_1 G_e = g_2G_e = g_3 G_e$ and $\tau \in S_3$ is the unique permutation for which $h_{\tau(1)} <_e h_{\tau(2)} <_e h_{\tau(3)}$ where 
$h_i = g^{-1} g_i$ for some (and hence any) $g \in g_1G_e$.  
\end{itemize}
It is routine to verify that $c'$ is a circular order on $G$. 
\end{proof}
A basic example is $E = S^1$ and $c$ the circular order determined by the usual orientation on $S^1$. In this case $\mbox{Aut}(E, c) = \mbox{Homeo}_+(S^1)$, and we assume below that $\mbox{Homeo}_+(S^1)$ is endowed with a circular order of the form described in the proof of Lemma \ref{lemma: aut is co}.

\subsection{Stir frying representations via recalibration} 
\label{subsec: pa flow action on the tree}

We saw in the previous subsection that if a group $G$ acts on a cyclically ordered $\mathbb R$-order tree $T$, there is an associated action of $G$ on the set of ends of $T$ which preserves the induced circular order. We will see in this subsection that this often leads to a circular order on $G$ via Lemma \ref{lemma: aut is co} and hence, if $G$ is countable, a dynamic realisation $\rho: G \to \G$. Altering the $G$-invariant cyclic ordering on $T$ by a process we call {\it recalibration}, yields a fundamentally new representation $G \to \G$. Here, we make this idea precise by applying it to the natural actions of the fundamental groups on the ends of the leaf spaces of stable laminations of pseudo-Anosov flows on universal covers.

As in \S \ref{sec: univ circs and locbcs hyp}, we assume that $\mathcal{M}$ is a closed, connected, oriented $3$-orbifold finitely covered by a manifold whose singular set is contained in an oriented link $B = B_1 \cup \cdots \cup B_m \subset |\mathcal{M}|$. Let $n_i \geq 1$ be the order of the isotropy of $\mathcal{M}$ along $B_i$. We denote a positively oriented meridional curve of $N(B_i)$ by $\beta_i$. 

Given a flow $\Phi$ on $|\mathcal{M}|$ which is well-adapted to the pair $(\mathcal{M}, B)$ (\S \ref{subsec: well adapted flows}), we constructed a commutative diagram of covering maps in the proof of Proposition \ref{prop: pa implies co orbifold}: 
\begin{center} 
\begin{tikzpicture}[scale=0.6]
\node at (8, 6) {$\widetilde{\mathcal{M}}$};
\node at (11.4, 3.5) {$W$};
\node at (8, 1) {$\mathcal{M}$}; 

\node at (7.65, 3.5) {$p$};  
\node at (10, 5.3) {$p'$};
\node at (10, 1.8) {$p_1$}; 
 
\draw [ ->] (8, 5.4) -- (8,1.5);
\draw [ ->] (8.5, 5.6) --(11,3.9); 
\draw [ ->] (10.9, 3.1) -- (8.6,1.4);
\end{tikzpicture}
\end{center} 
where $p$ and $p'$ are universal covers and $p_1$ is a finite degree regular cover from a manifold $W$ to $\mathcal{M}$. Let $\Phi'$ be the lift of the flow $\Phi$ to $W$ and $\widetilde \Phi$ its lift to $\widetilde{\mathcal{M}} \cong \mathbb{R}^3$. Both $\Phi'$ and $\widetilde{\Phi}$ are pseudo-Anosov.

The orbit space $\mathcal{O}$ of $\widetilde{\Phi}$ is homeomorphic to $\mathbb{R}^2$ (Theorem \ref{thm: orbit space}) and inherits an orientation induced from those on $\mathcal{M}$ and the flow. Let $\widetilde \Lambda_s$ on $\widetilde{\mathcal{M}} \cong \mathbb{R}^3$ be the pullback of the stable lamination $\Lambda_s'$ of $\Phi'$ (\S \ref{subsec: degeneracy loci}). It follows from our constructions that the image of $\widetilde \Lambda_s$ under the projection map $\pi: \widetilde{\mathcal{M}} \rightarrow \mathcal{O}$ is an essential lamination  $\bar \Lambda_s$ on $\mathcal{O}$ by lines. Then $\pi_1(\mathcal{M})$ acts on the leaf space $T(\bar \Lambda_s)$ and hence on its space of ends $\mathcal{E}(T(\bar \Lambda_s))$. Let
$$\varphi: \pi_1(\mathcal{M}) \to \mbox{Aut}(\mathcal{E}(T(\bar \Lambda_s)))$$
be the associated homomorphism. There is a $\pi_1(\mathcal{M})$-invariant cyclic ordering on $T(\bar \Lambda_s)$ induced by the inclusion $\bar \Lambda_s \subset \mathcal{O}$ so if $c_\Phi$ is the associated circular ordering on $\mathcal{E}(T(\bar \Lambda_s))$, the image of $\varphi$ is contained in the circularly ordered group $\mbox{Aut}(\mathcal{E}(T(\bar \Lambda_s)), c_\Phi)$. In particular, $\mbox{image}(\varphi)$ is a circularly ordered group.

Let $\Phi_0$ be the restriction of $\Phi$ to $|\mathcal{M}|\setminus B$. If for $i = 1, 2, \ldots, m$, we set $\Delta_i =  |\beta_i\cdot \delta_i(\Phi_0)|$, then each component $B_i$ of $B$ lifts to a countable union of flow lines in $\widetilde{\mathcal{M}}$ which corresponds to a $\pi_1(\mathcal{M})$-invariant subset $V_i$ of $T(\bar \Lambda_s)$, each point of which has valency 
$$d_i = n_i \Delta_i$$ 
For each $x \in V_i$, let $\sigma_{1}^x, \sigma_{2}^x, \ldots, \sigma_{d_i}^x$ be segments incident to $x$ such that the half-open intervals $\sigma_{i}^x \setminus \{x\}$ lie in different components of $T(\bar \Lambda_s)  \setminus \{x\}$. Assume, moreover, that they are indexed (mod $d_i$) with respect to the local circular  order at $x$ determined by the inclusion $\bar \Lambda_s \subset \mathcal{O}$. Then if $\beta_i^x \in \pi_1(\mathcal{M})$ is the conjugate of $\beta_i$ which leaves $x$ invariant, the proof of Proposition \ref{prop: pa implies co orbifold} shows that 
\begin{equation} 
\label{eqn: cyc order action} 
\beta_i^x \cdot \sigma_i^x = \sigma_{i+\Delta_i}^x
\end{equation}
Further, since each $\gamma \in \pi_1(\mathcal{M})$ acts as an orientation-preserving homeomorphism of $\mathcal{O}$, if $x \in V_i$ there is a $k(x, \gamma) \in \mathbb Z$ such that 
\begin{equation} 
\label{eqn: how transports} 
\gamma \cdot \sigma_j^x = \sigma_{j + k(x, \gamma)}^{\gamma \cdot x}
\end{equation}

The following lemma is needed for the proof of Theorem \ref{thm: pa implies co orbifold 2} below. 

\begin{lemma}
\label{lem: action faithful}
The kernel of $\varphi: \pi_1(\mathcal{M}) \to \mbox{Aut}(\mathcal{E}(T(\bar \Lambda_s)))$ is torsion-free. Moreover, its image is non-abelian and infinite if some $n_i \geq 3$. 
\end{lemma}

\begin{proof}
A non-trivial element $\gamma$ of finite order in $\pi_1(\mathcal{M})$ has fixed points in $\widetilde{\mathcal{M}} \cong \mathbb R^3$, so is a power of a conjugate of some $\beta_j$ where $n_j \geq 2$. Then $\gamma$ fixes a point $x$ of $T(\bar \Lambda_s)$ and by (\ref{eqn: cyc order action}) 
it acts non-trivially on $\mathcal{E}(T(\bar \Lambda_s))$. Hence $\gamma \not \in \mbox{kernel}(\varphi)$. 

Circularly ordered finite groups are cyclic, so to complete the proof we need only show that the image of $\varphi$ is non-abelian when some $n_i \geq 3$. 

By construction, $V_i$ is spread uniformly across $T(\bar \Lambda_s)$. In particular there are distinct points $x_0$ and $x_0' = \gamma \cdot x_0$, where $\gamma \in \pi_1(\mathcal{M})$. Then $ \beta_i^{x_0}$ stabilises $x_0$ while $ \beta_i^{x_0'} = \gamma  \beta_i^{x_0} \gamma^{-1}$ stabilises $x_0'$. If $e$ is an end of $T(\bar \Lambda_s)$ determined by a ray based at $x_0$ and passing through $x_0'$, the fact that $n_i \geq 3$ implies that $( \beta_i^{x_0}  \beta_i^{x_0'})(e) \ne ( \beta_i^{x_0'}  \beta_i^{x_0})(e)$. Hence $\varphi( \beta_i^{x_0}  \beta_i^{x_0'}) \ne \varphi( \beta_i^{x_0'}  \beta_i^{x_0})$, so the image of $\varphi$ is non-abelian.
\end{proof}

\begin{thm}
\label{thm: pa implies co orbifold 2}
Suppose that $\mathcal{M}$ is a closed, connected, oriented $3$-orbifold and $B = B_1 \cup \cdots \cup B_m$ is an oriented link in the underlying $3$-manifold $|\mathcal{M}|$ which contains the singular set of $\mathcal{M}$. Denote by $n_i \geq 1$ the order of the isotropy of $\mathcal{M}$ along $B_i$ and suppose that $\Phi$ is a flow on $|\mathcal{M}|$ that is well-adapted to $(\mathcal{M}, B)$. Then given integers $a_i$ coprime with $n_i $, there is an induced action of $\pi_1(\mathcal{M})$ on a cyclically ordered $\mathbb R$-order tree and an associated faithful representation  
$$\rho: \pi_1(\mathcal{M}) \to \mbox{{\rm Homeo}}_+(S^1)$$
where $\rho( \beta_i)$ is conjugate to rotation by $2 \pi a_i/n_i$ for each $i$. 
\end{thm}

We note that the theorem allows for $n_i$ to be $1$ and that if $n_i=1$ for all $i$, then $\mathcal{M}$ is a manifold.

The following lemma will be used in the proof of the theorem. 

\begin{lemma}
\label{lemma: coprime with deltai}
Let $a_i$ and $n_i$ be as in the theorem and let $\Phi_0$ be the pseudo-Anosov flow obtained by restricting $\Phi$ to $C(B)$.  If $\Delta_i = |\delta_i(\Phi_0) \cdot \beta_i| \geq 1$ then there exists $a_i' \equiv a_i \mbox{ {\rm (mod $n_i$)}}$ such that $\gcd(a_i', \Delta_i) = 1$. 
\end{lemma}

\begin{proof}
Factor $\Delta_i$ as $c_i e_i$ where $\gcd(e_i, n_i) = 1$ and each prime which divides $c_i$ also divides $n_i$. Next choose an integer $k$ so that $kn_i \equiv 1 - a_i$ (mod $e_i$), so we have $a_i' = a_i + kn_i \equiv 1$ (mod $e_i$). Then $a_i' \equiv a_i$ (mod $n_i$) and is coprime with $e_i, n_i$ and, a fortiori, with $c_i$. Thus it is coprime with $\Delta_i = c_i e_i$. 
\end{proof}

\begin{proof}[Proof of Theorem \ref{thm: pa implies co orbifold 2}]
Let $\Phi_0$ denote the restriction of $\Phi$ to $C(B)$, which, by hypothesis, is pseudo-Anosov. By Lemma \ref{lemma: basic top consequences of cpaf},  $\mathcal{M}$ contains no teardrops or spindles, so is finitely covered by a manifold \cite{BLP05}.

If $n_i = 1$ or $2$ for all $i=1, \cdots, m$, we take $\rho$ to be the asymptotic circle action $\rho$ in Proposition \ref{prop: pa implies co orbifold}. We remark that we cannot always use the action of $\pi_1(\mathcal{M})$ on $\mathcal{E}(T(\bar \Lambda_s))$ in this case as it is possible that $T(\bar \Lambda_s)$ is a line and hence $\mathcal{E}(T(\bar \Lambda_s))$ only contains two points.

Next we assume that some $n_i \geq 3$. Lemma \ref{lem: action faithful} then shows that the kernel of $\varphi: \pi_1(\mathcal{M}) \to \mbox{Aut}(\mathcal{E}(T(\bar \Lambda_s)))$ is torsion free and its image is infinite and non-abelian. 

Our strategy involves altering the circular orderings on $T(\bar \Lambda_s)$ equivariantly over each $V_i$ with $n_i \geq 3$. (When $n_i\leq 2$, $a_i$ is uniquely determined (mod $n_i$) and no such alterations are needed.) The reader will find an illustrative example of the idea of the argument in Figure \ref{fig: circular order ends} below. 

Recall that $V_i$ is the $\pi_1(\mathcal{M})$-invariant subset of $T(\bar \Lambda_s)$ corresponding to $B_i$ and that each point of $V_i$ has valency $d_i = n_i \Delta_i$. Without loss of generality we can suppose that $\gcd(a_i, d_i) = 1$, by Lemma \ref{lemma: coprime with deltai}.   

For each $i$ with $n_i \geq 3$, fix an integer $b_i$ such that $a_i b_i \equiv 1$ (mod $d_i$) and replace the circular ordering at $x \in V_i$ by the one determined by the listing 
\begin{eqnarray} 
\label{eqn: new local order} 
\sigma_1^x, \sigma_{b_i + 1}^x,  \sigma_{2b_i + 1}^x, \sigma_{3b_i + 1}^x, \ldots, \sigma_{(d_i-1)b_i + 1}^x 
\end{eqnarray}
where the indices take values in $\{1, \cdots, d_i\}$ (mod $d_i$).  
We claim that this new circular ordering on the segments incident to the points of $V_i$ is invariant under the action of $\pi_1(\mathcal{M})$. To see  this, first use (\ref{eqn: cyc order action}) and (\ref{eqn: how transports}) to reduce the verification to showing that the new local order is preserved under the action $t: \{\sigma_i^x\} \rightarrow \{\sigma_i^x\}$ with $t(\sigma_i^x) = \sigma_{i+1}^x$. Next note that since $a_i b_i \equiv 1$ (mod $d_i$), for any $k \in \{0, 1, \cdots, d_i - 1\}$, we have 
\begin{equation}
\label{equ: actions on new local order}
    \sigma_{(kb_i + 1) + 1}^x =  \sigma_{(kb_i + 1) + a_ib_i}^x = \sigma_{(k+a_i)b_i + 1}^x 
\end{equation}
and therefore the sequence 
\begin{displaymath}
    t(\sigma_{1}^x), t(\sigma_{b_i + 1}^x),  t(\sigma_{2b_i+1}^x),  \ldots, t(\sigma_{(d_i-1)b_i +1}^x)
\end{displaymath}
is identical to the sequence in (\ref{eqn: new local order}) up to a cyclic permutation.

We {\em recalibrate} $T(\bar \Lambda_s)$ by changing the local order at $V_i$ as in (\ref{eqn: new local order}) to obtain a new $\pi_1(\mathcal{M})$-invariant cyclic ordering on $T(\bar \Lambda_s)$ and therefore a new $\pi_1(\mathcal{M})$-invariant circular order $c$ on $\mathcal{E}(T(\bar \Lambda_s))$. As such, the image of $\varphi: \pi_1(\mathcal{M}) \to \mbox{Aut}(\mathcal{E}(T(\bar \Lambda_s)))$ lies in $\mbox{Aut}(\mathcal{E}(T(\bar \Lambda_s)), c)$. Since $\mbox{kernel}(\varphi)$ is a torsion-free infinite index subgroup of $\pi_1(\mathcal{M})$, it is the fundamental group of an irreducible, orientable, non-compact $3$-manifold. Hence it is either trivial or left-orderable (cf. proof of \cite[Theorem 1.1]{BRW05}). Then the exact sequence
$$1 \to \mbox{kernel}(\varphi) \to \pi_1(\mathcal{M}) \xrightarrow{\; \varphi \;} \mbox{Aut}(\mathcal{E}(T(\bar \Lambda_s)), c)$$
determines a circular order $c_0$ on $\pi_1(\mathcal{M})$ (\cite[Lemma 2.2.12]{Cal04}) whose dynamic realisation is an order-preserving injection $\rho_0: \pi_1(\mathcal{M}) \to \mbox{Homeo}_+(S^1)$ (\cite[Lemma 2.2.10]{Cal04}). 

We show that $\rho_0(\beta_i^x)$ is conjugate to a rotation by $2 \pi a_i/n_i$. To do this, we must unpack the construction of the circular order $c_0$ on $\pi_1(\mathcal{M})$ and its dynamic realisation $\rho_0$.

The set of ends $\mathcal{E}(T(\bar \Lambda_s))$ can be decomposed as the disjoint union $\bigsqcup_{j=1}^{d_j} \mathcal{E}_j^x$, where $\mathcal{E}_j^x$ is the set of ends determined by infinite geodesic spines based at $x$ with initial segment $\sigma_j^x$. Then (\ref{eqn: cyc order action}) shows that $\beta_i^x \cdot \mathcal{E}_j^x = \mathcal{E}_{j+\Delta_i}^x$

Choose $e_j \in \mathcal{E}_j^x$ for $j = 1, 2, \ldots, d_i$ so that the set $\{e_j\}$ is preserved under the action of $\beta_i^x$. By construction, the circular order on the set $\{e_j\}$ corresponds to the circular order on the set of segments $\{\sigma_j^x\}$, which is given by (\ref{eqn: new local order}). See Figure \ref{fig: circular order ends}. 

\begin{figure}[ht]
    \centering    
    \includegraphics[scale=0.7]{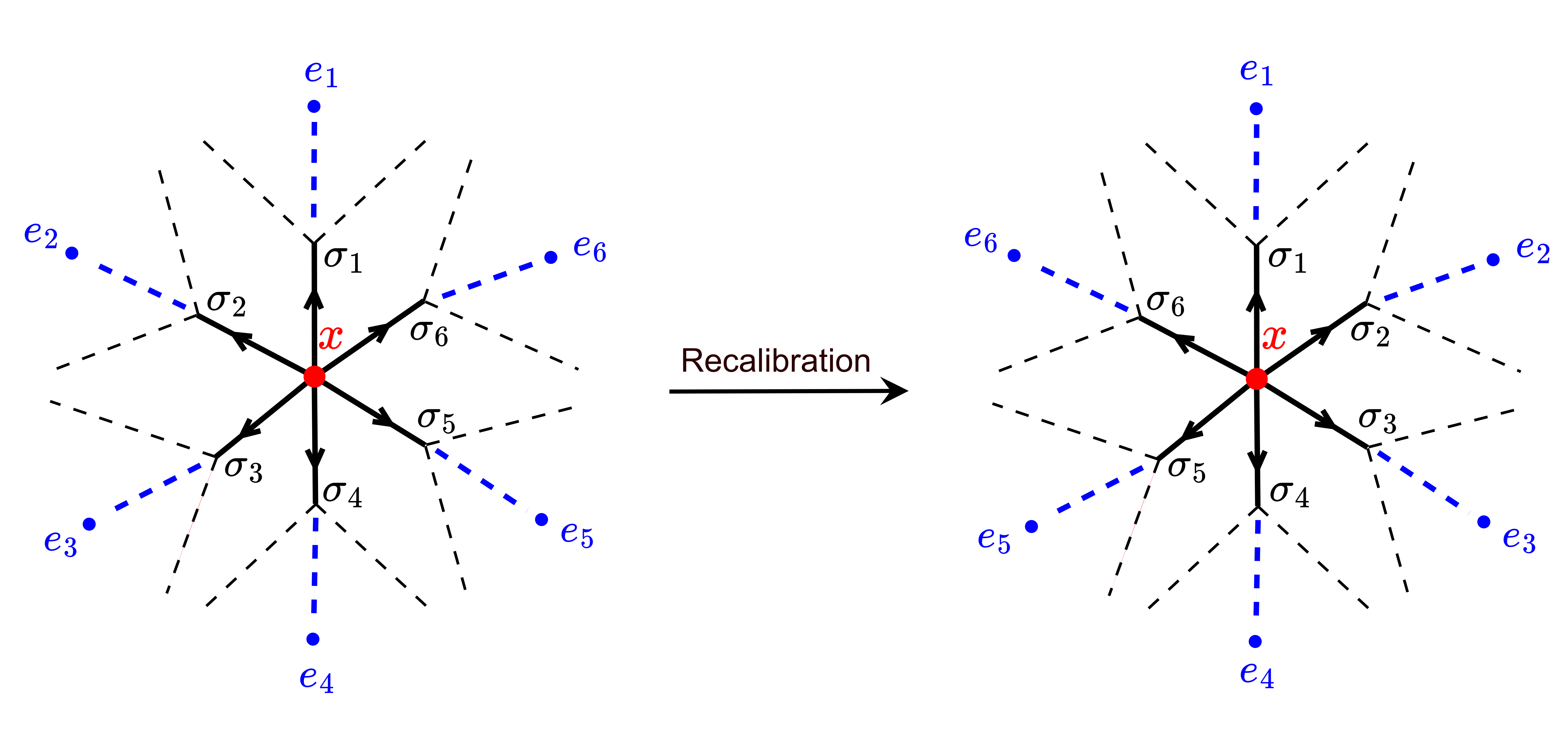}
    \caption{In this figure we simplify the notation by dropping the superscript $x$ from $\sigma_j^x$. We suppose that $n_i=3$, $\Delta_i = |\delta_{i}(\Phi_0) \cdot \beta_i| = 2$, so $d_i = 2\times 3 = 6$, and $a_i = 2$. Theorem \ref{thm: pa implies co orbifold 2} claims that we can find a representation $\rho$ so that $\rho(\beta_i)$ is conjugate to rotation by $4\pi/3$. Since $(a_i, d_i)\neq 1$, to make the construction work, we must replace $a_i$ by $5$, say, which equals $2$ (mod $n_i$), and is coprime with $d_i = 6$ (cf. Lemma \ref{lemma: coprime with deltai}). Choose $b_i = -1$. Now we can produce the new local order as described in (\ref{eqn: new local order}), which is shown in the figure to the right. Note that since $\Delta_i = 2$, $\beta_i$ sends $\sigma_i^x$ to $\sigma_{i+2}^x$ so its action under the natural order has rotation number $1/3$ (mod $\mathbb{Z}$), while its rotation number under the new order is $2/3$ (mod $\mathbb{Z}$).}
    \label{fig: circular order ends}
\end{figure}

Given any $k\in \{0,\ldots, d_i-1\}$, by (\ref{equ: actions on new local order}), we have
\begin{equation}
    \label{equ: mu action on ends}
     \beta_i^x\cdot e_{1 + kb_i} = e_{1 + kb_i + \Delta_i} = e_{1 + (k + a_i \Delta_i)b_i} 
\end{equation}
If we use $e_1$ to define the circular order on $\mbox{Aut}(\mathcal{E}(\bar \Lambda_s), c)$ (cf. the proof of Lemma \ref{lemma: aut is co}), then the circular order on the cyclic group $\langle  \beta_i^x \rangle$ is determined by the induced circular order on the set $\langle  \beta_i^x\rangle\cdot e_1 \subset \mathcal{E}(\bar \Lambda_s)$.

Choose an order-preserving embedding $\iota: \pi_1(\mathcal{M}) \to S^1$ such that the action of $\pi_1(\mathcal{M})$ on $\iota(\pi_1(\mathcal{M}))$ given by $\gamma' \cdot  \iota(\gamma) = \iota(\gamma' \gamma)$ for any $\gamma, \gamma'\in\pi_1(\mathcal{M})$, extends to an action of $\pi_1(\mathcal{M})$ on $S^1$ by orientation-preserving homeomorphisms. Assuming that $\iota$ is chosen with enough care (cf. \cite[Lemma 2.2.10]{Cal04}), we obtain the dynamic realisation of $c_0$, which is a faithful representation $\rho_0: \pi_1(\mathcal{M}) \to \mbox{Homeo}_+(S^1)$.

Then from (\ref{equ: mu action on ends}) and the construction of the dynamic realization, it follows that the rotation number of $\rho_0( \beta_i^x)$ in ${\rm Homeo}_+(S^1)$ is $\frac{a_i \Delta_i}{d_i} = \frac{a_i}{n_i}$.  Since $ \beta_i^x$ has finite order, it must be conjugate to a rotation by $\frac{2\pi a_i }{n_i}$, which was to be proved.
\end{proof}

\begin{thm} 
\label{thm: reps result intro} 
Let $L = K_1 \cup \cdots \cup K_m$ be a link in an orientable $3$-manifold $W$ whose complement admits a pseudo-Anosov flow $\Phi_0$. For each $i$, fix an  oriented essential simple closed curve $\alpha_i$ on $T_i$ and an integer $n_i \geq 1$ so that $n_i|\alpha_i \cdot \delta_i(\Phi_0)| \geq 2$. Then for any integer $a_i$ coprime with $n_i$, there is a homomorphism $\rho: \pi_1(X(L)) \to \mbox{Homeo}_+(S^1)$ with non-cyclic image such that $\rho(\alpha_i)$ is conjugate to rotation by $2\pi a_i/n_i$ for each $i$.
\end{thm}

\begin{proof}
Set $\mathcal{M} = X(L)(\alpha_*; n_*)$. Let $B = B_1 \cup \cdots \cup B_m$ be the link in $|\mathcal{M}| = X(L)(\alpha_*)$ corresponding to the cores of the $\alpha_i$-surgery solid tori. Theorem \ref{thm: reps result intro} now follows from Theorem \ref{thm: pa implies co orbifold 2} as Corollary \ref{cor: branched covers of pseudo-Anosov links} followed from Proposition \ref{prop: pa implies co orbifold}. 
\end{proof}

\section{Euler class of representations and left-orderable cyclic branched covers}
\label{sec: reps and locbcs}

The goal of this section is to prove the following theorem.  Our method is to show that the Euler classes, defined below, of the representations constructed in the previous sections vanish on certain co-cyclic subgroups.

\begin{thm} 
\label{thm: locbc intro}
Let $L = K_1 \cup \cdots \cup K_m$ be a prime link in an integer homology $3$-sphere whose exterior is irreducible. Suppose that $\rho: \pi_1(X(L)) \to \mbox{Homeo}_+(S^1)$ is a representation with non-cyclic image such that $\rho(\mu_i)$ is conjugate to rotation by $2\pi a_i/n$ for some $a_i, n\in \mathbb{Z}$, where $n \geq 2$. If the induced homomorphism $\psi: \pi_1(X(L)) \to \mathbb Z/n$ which sends $\mu_i$ to $a_i$ $(${\rm{mod}} $n)$ is an epimorphism, then $\pi_1(\Sigma_\psi(L))$ is left-orderable. 
\end{thm}

\begin{remark}
We allow the possibility that $\rho(\mu_i) = \mbox{id}_{S^1}$ in Theorem \ref{thm: locbc intro} (i.e. $a_i \equiv 0$ (mod $n$)), which leads to a slightly more general type of cyclic branched cover $\Sigma_\psi(L)$ than considered elsewhere in the paper. Similarly we allow $n_i$ to be $1$ in Lemma \ref{lemma: h^2} below. 
\end{remark}

\subsection{Euler classes} 
\label{subsec: Euler classes}
The set of central extensions of a group $G$ by $\mathbb Z$ is naturally identified with $H^2(G)$ in such a way that the direct product $G \times \mathbb Z$ corresponds to $0$. More precisely, each inhomogeneous $2$-cocycle $\xi$ on $G$ normalised to take the value $0$ on $(1,1)$ determines a central extension 
$$1 \to \mathbb Z \to \widetilde G_\xi \xrightarrow{\; \varphi \;} G \to 1,$$
where $\widetilde G_\xi = G \times \mathbb Z$ as a set, $\varphi$ is the projection, and multiplication in $\widetilde G_\xi$ is defined by $(g, a) \cdot (h, b) = (gh, a + b + \xi(g, h))$. Altering the cocycle by a coboundary yields an equivalent extension. Conversely, given a central extension $1 \to \mathbb Z \to \widetilde G \xrightarrow{\; \varphi \;} G \to 1$ and transversal $s: (G, 1) \to (\widetilde G, 1)$ to $\varphi$, the function $\xi: G^2 \to  \mbox{kernel}(\varphi) = \mathbb Z, (g, h) \mapsto s(gh)^{-1}s(g)s(h)$ is a normalised $2$-cocycle whose class in $H^2(G)$ is independent of the choice of $s$.  

The class $e \in H^2(G)$ of an extension $1 \to \mathbb Z \to \widetilde G \xrightarrow{\; \varphi \;} G \to 1$ is called its {\it Euler class}. 

There is a universal covering homomorphism $\varphi: \mbox{{\rm Homeo}}_{\mathbb Z}(\mathbb R) \to \mbox{{\rm Homeo}}_+(S^1)$, where $\mbox{{\rm Homeo}}_{\mathbb Z}(\mathbb R)$ is the group of homeomorphisms of the real line which commute with translation by $1$. The kernel of $\varphi$ is the group of integer translations and is central in $\mbox{{\rm Homeo}}_{\mathbb Z}(\mathbb R)$, so there is a central extension 
$$1 \to \mathbb Z \to \mbox{{\rm Homeo}}_{\mathbb Z}(\mathbb R) \xrightarrow{\; \varphi \;} \mbox{{\rm Homeo}}_+(S^1) \to 1$$
Given a representation $\rho: G \to \mbox{Homeo}_+(S^1)$, define $\widetilde G_{\rho}$ to be the subgroup $\{(g, f) \; | \; \rho(g) = \varphi(f)\}$ of the direct product $G \times \mbox{{\rm Homeo}}_{\mathbb Z}(\mathbb R)$ and note that the projections $\varphi_\rho: \widetilde G_{\rho} \to G$ and $\tilde \rho: \widetilde G_{\rho} \to \mbox{{\rm Homeo}}_{\mathbb Z}(\mathbb R)$ give rise to a commutative diagram of central extensions 

\begin{center} 
\begin{tikzpicture}[scale=0.8]
 
\node at (7, 4.5) {$\mathbb Z$};
\node at (10.85, 4.5) {$\widetilde G_{\rho}$};
\node at (15.1, 4.5) {$G$};
 
\node at (7, 1.5) {$\mathbb Z$}; 
\node at (11, 1.5) {$\mbox{{\rm Homeo}}_{\mathbb Z}(\mathbb R)$};
\node at (15.1, 1.5) {$\mbox{{\rm Homeo}}_+(S^1)$};

\draw [ >->] (7.4, 4.5) --(10.35,4.5); 
\draw [ ->>] (11.5, 4.5) --(14.7,4.5); 
 
\draw [ ->] (7, 4) -- (7,2);
\draw [ ->] (10.8, 4) -- (10.8,2);
\draw [ ->] (15.1, 4) -- (15.1,2); 
 
\draw [ >->] (7.4, 1.5) --(9.6,1.5); 
\draw [ ->>] (12.35, 1.5) --(13.55,1.5); 

\node at (6.6, 3) {${\rm {\tiny id}}$};
\node at (10.4, 3) {$\tilde \rho$};
\node at (14.6, 3) {$\rho$};

\node at (13, 5) {$\varphi_\rho$};
\node at (13, 2) {$\varphi$};
 
\end{tikzpicture}
\end{center} 
The {\it Euler class} $e(\rho) \in H^2(G)$ of $\rho$ is defined to be the Euler class of the extension
$1 \to \mathbb Z \to  \widetilde G_{\rho} \xrightarrow{\; \varphi_\rho \;} G \to 1$  
and is therefore the obstruction to $\varphi_\rho$ admitting a splitting homomorphism. Equivalently, it is the obstruction to lifting $\rho$ to a representation $G \to \mbox{Homeo}_{\mathbb Z}(\mathbb R)$.

\subsection{Proof of Theorem \ref{thm: locbc intro}}
\label{subsec: reps and locbcs}
We begin with a lemma. 

\begin{lemma}
\label{lemma: h^2}
Suppose that $L$ is an oriented link in an oriented integer homology $3$-sphere $W$ with components $K_1, K_2, \ldots, K_m$ and $\psi: \pi_1(X(L)) \to \mathbb Z/n$ is an epimorphism with associated $n$-fold cyclic branched cover $(\Sigma_\psi(L), \hat L) \to (W, L)$. If $n_i$ is the order of $\psi(\mu_i)$ in $\mathbb Z/n$ and $\Sigma_\psi(L)$ is an irreducible rational homology $3$-sphere with infinite fundamental group, then $H^2(\pi_1(X(L)(\mu_*; n_*))) \cong \oplus_{i=1}^m \mathbb Z/n_i$, where $\mu_* = (\mu_1,\mu_2, \ldots, \mu_m)$ and $n_* = (n_1, n_2, \ldots, n_m)$.
\end{lemma}

\begin{proof}
Set $G = \pi_1(X(L)(\mu_*; n_*))$. Our hypotheses imply that $\psi$ factors through a homomorphism $\bar \psi: G \to \mathbb Z/n$ where, by construction, the cover of $X(L)(\mu_*; n_*)$ corresponding to the kernel of $\bar \psi$ is $\Sigma_\psi(L)$. Let $Z$ be a $K(G, 1)$ so that $H_*(G) \cong H_*(K)$ and $H^*(G) \cong H^*(Z)$. Since $W$ is an integer homology $3$-sphere, it is easy to verify that $H_1(Z) \cong \oplus_{i=1}^m \mathbb Z/n_i$

Let $\widetilde Z \to Z$ be the $n$-fold cyclic cover where $\pi_1(\widetilde Z) = \pi_1(\Sigma_\psi(L)) = \mbox{kernel}(\bar \psi)$, so that $\widetilde Z = K(\pi_1(\Sigma_\psi(L)), 1)$. Since $\Sigma_\psi(L)$ is irreducible with infinite fundamental group, it is aspherical. Hence $\Sigma_\psi(L) \simeq \widetilde Z$ and therefore as $\Sigma_\psi(L)$ is a rational homology $3$-sphere, $H_2(\widetilde Z) = H_2(\Sigma_\psi(L)) = 0$. A transfer argument then shows that 
$H_2(Z; \mathbb Q) = 0$. Now apply universal coefficients to $Z$ to deduce that
$$H^2(G) \cong H^2(Z) \cong \mbox{Ext}(H_1(Z),\mathbb Z) \oplus \mbox{Hom}(H_2(Z), \mathbb Z) = \mbox{Ext}(\oplus_{i=1}^m \mathbb Z/n_i,\mathbb Z) \cong \oplus_{i=1}^m \mathbb Z/n_i$$
\end{proof}

\begin{proof}[Proof of Theorem \ref{thm: locbc intro}] 
Let  $n_i$ be the order of $\psi(\mu_i)$ in $\mathbb Z/n$, $n_* = (n_1, n_2, \ldots, n_m)$, $\mu_* = (\mu_1,\mu_2, \ldots, \mu_m)$ and $G = \pi_1(X(L)(\mu_*; n_*))$. Our hypotheses imply that $\psi$ and $\rho$ factor through homomorphisms $\bar \psi: G \to \mathbb Z/n$ and $\bar \rho: G \to \mbox{{\rm Homeo}}_+(S^1)$ where, by construction, the cover of $X(L)(\mu_*; n_*)$ corresponding to $\bar \psi$ is $\Sigma_\psi(L)$. Hence $\pi_1(\Sigma_\psi(L)) = \mbox{ker}(\bar \psi)$. Our hypotheses also imply that $\Sigma_\psi(L)$ is irreducible (cf. \cite[Proposition 10.2]{BGH21}), so it will have a left-orderable fundamental group if it has a positive first Betti number (\cite[Theorem 1.1]{BRW05}). Assume, then, that $\Sigma_\psi(L)$ is a rational homology $3$-sphere and therefore $H^2(G) \cong \oplus_{i=1}^m \mathbb Z/n_i$ by Lemma \ref{lemma: h^2}.  

Since the only finite subgroups of $\mbox{{\rm Homeo}}_+(S^1)$ are conjugate into $SO(2)$, hence cyclic, our hypotheses also imply that the image of $\bar \rho$ is infinite. Then the restriction of $\bar \rho$ to $\pi_1(\Sigma_\psi(L))$ also has infinite image. In particular, $\Sigma_\psi(L)$ has an infinite fundamental group. 

Let $Y$ be a CW-complex obtained by attaching $2$-cells $D_1, D_2, \ldots, D_m$ to $X(L)$ with attaching maps $\partial D_i \to X(L)$ which wrap $n_i$ times around a positively oriented meridian of $K_i$, so that $\pi_1(Y) = G$. Next construct an Eilenberg-MacLane space $Z = K(G, 1)$ by attaching cells of dimension $3$ or more to $Y$. We can identify $H^2(G)$ with $H^2(Z)$. 

In what follows we use $C_*(Z)$ and $C^*(Z)$ to denote the cellular chain and cochain complexes of $Z$ over $\mathbb Z$. A cohomology class represented by a cocycle $\zeta$ will be denoted by $[\zeta]$.  

For each $1 \leq i \leq m$, let $f_i \in \mbox{Homeo}_\mathbb Z(\mathbb R)$ be a lift of $\bar \rho(\mu_i)$ which is conjugate to a translation by $a_i/n$. In \S 2 of \cite{Mil58}, Milnor constructed a cellular $2$-cocycle $\omega$ of $Z$ representing $e(\bar \rho)$ whose value on the $2$-cell $D_i$ is the negative of the translation number of $f_i^{n_i} = \mbox{sh}(n_ia_i/n)$ (\cite[Lemma 2]{Mil58}). That is, $\omega(D_i) = -n_ia_i/n$, where we note that as $n_i$ is the order of $\psi(\mu_i) \equiv a_i$ in $\mathbb Z/n$, $n_ia_i/n \in \mathbb Z$. 

Since $H^2(Z) \cong H^2(G) \cong \oplus_{i=1}^m \mathbb Z/n_i$ where each $n_i$ divides $n$, $n \omega$ represents zero in $H^2(Z)$. Hence there is a $1$-cochain  $\eta \in C^1(Z)$ for which $\delta \eta = n\omega$ and so reducing $\eta$ (mod $n$) we obtain a (mod $n$) 1-cocycle $\bar \eta$ representing an element $[\bar \eta] \in H^1(Z; \mathbb{Z}/n) = \mbox{Hom}(H_1(Z), \mathbb Z/n)$. Let $f$ be the composition 
$$\pi_1(Z) \rightarrow H_1(Z) \xrightarrow{\; [\bar\eta] \;} \mathbb{Z}/n$$ 
and note that if $p: \widetilde{Z} \rightarrow Z$ is the covering corresponding to $\mbox{kernel}(f)$, then $\widetilde{Z} = K(\mbox{kernel}(f), 1)$. We claim that $\widetilde Z$ is homotopy equivalent to $\Sigma_{\psi}(L)$. Equivalently, $\pi_1(\widetilde{Z}) \cong \pi_1(\Sigma_{\psi}(L))$. 

To see this, note that the boundary of the $2$-chain $D_i$ is the $1$-cycle $n_i  \mu_i$ in $C_1(Z)$, so 
\begin{equation}
\eta(\mu_i) = (1/n_i) \eta (\partial (D_i)) = (1/n_i)(\delta \eta) (D_i) = (n/n_i)\omega (D_i) = (n/n_i)(-n_ia_i/n) = -a_i\nonumber 
\end{equation}
Then thinking of $\mu_i$ as an element of $\pi_1(Z) = G$ we have 
$$f(\mu_i) = \bar \eta(\mu_i) \equiv -a_i \mbox{ (mod $n$) } = - \bar \psi(\mu_i)$$
It follows that $f = -\bar \psi$ and therefore $\pi_1(\widetilde{Z}) = \mbox{kernel}(f)=  \mbox{kernel}(\bar \psi) = \pi_1(\Sigma_\psi(L))$.  

To complete the proof, let $\hat \rho = \bar \rho|_{\pi_1(\Sigma_\psi(L))}$ and observe that  
$$e(\hat \rho) = p^*(e(\bar \rho)) = p^*([\omega]) \in H^2(\widetilde{Z})$$
If we can show that $p^*([\omega]) = 0$, then $\hat \rho$ lifts to a (non-trivial) representation $\tilde \rho: \pi_1(\Sigma_\psi(L)) \to \mbox{Homeo}_\mathbb Z(\mathbb R) \leq \mbox{Homeo}_+(\mathbb R)$, so as $\Sigma_\psi(L)$ is orientable and irreducible, $\pi_1(\Sigma_\psi(L))$ is left-orderable  (\cite[Theorem 1.1]{BRW05}) and we are done. 

To prove that $p^*([\omega]) = 0$, note that the identity $\delta \eta = n \omega$ implies that the Bockstein homomorphism $H^1(Z; \mathbb Z/n) \xrightarrow{\; \beta \;} H^2(Z)$ of the coefficient sequence $0 \to \mathbb Z \xrightarrow{\; n \;} \mathbb Z \to \mathbb Z/n \to 0$ sends $[\bar \eta]$ to $[\omega]$. Then
$$p^*([\omega]) = p^*(\beta([\bar \eta])) = \beta(p^*([\bar \eta]))$$
by the naturality of $\beta$. On the other hand, $p^*([\bar \eta]) \in H^1(\widetilde Z; \mathbb Z/n)$ corresponds to $[\bar \eta] \circ p_* \in \mbox{Hom}(H_1(\widetilde Z), \mathbb Z/n)$, so is zero by the definition of $\widetilde Z$. Thus $p^*([\omega]) = 0$. 
\end{proof}

\section{Applications and Examples}
\label{sec: examples}

\subsection{Order detects meridional slopes}
\label{subsec: app slope detection}

There are three types of slope detection: order-detection, foliation-detection and non-$L$-space detection. For simplicity, we consider slope detection for knot manifolds, i.e. compact, connected, orientable, irreducible $3$-manifolds  with torus boundary that are not homeomorphic to $S^1\times D^2$. It was shown in \cite[Theorem 1.3]{BC2} that given two knot manifolds $M_1$ and $M_2$ and a homeomorphism $f: \partial M_1 \rightarrow \partial M_2$, if $f$ maps an {\em order-detected} slope on $\partial M_1$ to an order-detected slope on $\partial M_2$, then the fundamental group of $W = M_1\cup_f M_2$ is left-orderable. Similar results hold for foliation-detection and non-$L$-space detection, which are proved in \cite[Theorem 5.2]{BGH21} and \cite[Theorem 1.14]{HRW1} respectively.

We refer the reader to \cite[\S 7.2]{BC2} (see also \cite[\S 6]{BGH21}) for the formal definition of order-detection. We remark that order-detection is called LO-detection in \cite{BGH21}. The following proposition gives a sufficient representation-theoretic condition for a slope to be order-detected.  

\begin{prop}[Proposition 6.9 in \cite{BGH21}]
\label{prop: order detection}
Suppose that $M$ is a knot manifold and $\alpha$ is an oriented essential simple closed curve on $\partial M$. If $\rho: \pi_1(M) \rightarrow \text{Homeo}_+(\mathbb{R})$ is a homomorphism such that $\rho(\alpha)$ has a fixed point but not $\rho(\pi_1(\partial M))$, then the slope given by $\alpha$ is order-detected.
\end{prop}

\begin{thm}
\label{thm: meridian detection}
Let $K$ be a hyperbolic knot in an integer homology sphere $W$, such that the complement of $K$ admits a pseudo-Anosov flow whose degeneracy locus is meridional. Then the knot meridian is order-detected. 
\end{thm}

\begin{proof}
Let $\mu, \lambda$ be simple closed meridional and longitudinal curves of $K$ which are arbitrarily oriented. By assumption, there is a pseudo-Anosov flow $\Phi_0$ on $C(K)$, such that up to flipping the orientation on $\mu$, satisfies $\delta(\Phi_0) = c\mu$ for some $c\geq 1$. 

Let $\alpha$ be an oriented simple closed curve representing the homology class $\mu + \lambda$ in $H_1(\partial X(K))$. For each $n \geq 2$ we consider the orbifold $\mathcal{M} = X(K)(\alpha; n)$ and let $B$ be the core of the filling solid torus. Since $n |\delta(\Phi_0) \cdot \alpha| = nc \geq  2$, the flow $\Phi_0$ extends to a flow $\Phi$ on the underlying manifold $|\mathcal{M}|$ which is well-adapted to the pair $(\mathcal{M}, B)$, at least after we orient $B$ with the induced orientation from the flow $\Phi$. 

Then by Theorem \ref{thm: pa implies co orbifold 2}, there exists a representation $\rho: \pi_1(\mathcal{M})\rightarrow \G$ such that $\rho(\alpha)$ is conjugate to a rotation by $\frac{2\pi}{n}$.  For the rest of the proof we use the notation established in the proof of Theorem \ref{thm: pa implies co orbifold 2} (\S \ref{subsec: pa flow action on the tree}). 

Recall that $\Phi$ lifts to a pseudo-Anosov flow on the universal cover $\widetilde {\mathcal{M}}$. Let $\widetilde{\Lambda}_s$ be its stable lamination, which projects to an essential lamination on the leaf space $\mathcal{O}$ of the pullback flow on the universal cover $\widetilde{\mathcal{M}}$, and $T(\bar\Lambda_s)$ denote the leaf space of $\bar \Lambda_s$ which is a cyclically ordered $\mathbb{R}$-order tree. In the case that $n>2$, the representation $\rho$ is from the action of $\pi_1(\mathcal{M})$ on $\mathcal{E}(T(\bar \Lambda_s))$, the set of ends of $T(\bar\Lambda_s)$.

The knot $K$ lifts to a countable union of flow lines in the universal cover $\widetilde{\mathcal{M}}$, which corresponds to a $\pi_1(\mathcal{M})$-invariant subset $V$ of $T(\bar\Lambda_s)$, each point of which has valency $nc$. For each $x\in V$, we let $\sigma_1^x, \cdots, \sigma_{nc}^x$ be the distinct segments of $T(\bar\Lambda_s)$ incident to $x$, which are indexed (mod $nc$) with respect to the local circular order at $x$. 

Let $\mu^x$ be the conjugate of $\mu$ in $\pi_1(\mathcal{M})$ which fixes $x$. The rotation number of $\rho(\mu^x)$ is determined by how $\mu^x$ acts on $\sigma_i^x$, which is dynamically identical with the action of $\mu^x$ on the cusps of the complementary region of $\widetilde{\Lambda}_s$ that is fixed by it. Since $\mu$ is parallel to a component of the degeneracy locus, it follows that $\mu^x$ fixes each $\sigma_i^x$ and therefore, the rotation number of $\rho(\mu^x)$ is zero. By 
Lemma \ref{lemma: rot and trans}, the rotation number of $\rho(\mu)$ is also zero.

Let $\rho'$ be the composition $\pi_1(X(K))\rightarrow \pi_1(\mathcal{M}) \rightarrow \G$. Since $H^2(X(K)) = 0$ and $H_1(X(K))$ is generated by the class represented by $\mu$, there exists a lift of $\rho'$, denoted by $\tilde{\rho}': \pi_1(X(K)) \rightarrow \tG$, such that the translation number of $\tilde{\rho}'(\mu) = 0$. However, the translation number of $\tilde{\rho}'(\alpha)$ equals $\frac{1}{n} + k$ for some $k\in\mathbb{Z}$ and hence is nonzero. It follows that $\tilde{\rho}'(\alpha)$ acts without fixed points on $\mathbb R$. By Lemma \ref{lemma: rot and trans} and Proposition \ref{prop: order detection}, we have that $\mu$ is order-detected.  
\end{proof}

\subsection{Orderability of branched covers and link orientations}
\label{subsec: lo and link orientation}

We begin with a proof of Theorem \ref{thm: lo cbcs intro}.

\begin{proof}[Proof of Theorem \ref{thm: lo cbcs intro}]
Let $\psi(\mu_i) = a_i$ (mod $n$) in $\mathbb{Z}/n$, $n_i = n/\gcd(a_i, n)$ and $a_i' = a_i/\gcd(a_i, n)$. Then $a_i'$ and $n_i$ are coprime, and $n_i$ is the order of $a_i$ in $\mathbb{Z}/n$ for each $i$. By assumption, $n_i \geq 2$, so as the degeneracy loci $\delta_i(\Phi_0)$ are non-meridional, we have $n_i|\mu_i\cdot \delta_i(\Phi_0)| \geq 2$. Then by Theorem \ref{thm: reps result intro}, there exists a homeomorphism $\rho: \pi_1(X(L)) \to \mbox{Homeo}_+(S^1)$ with non-cyclic image such that $\rho(\mu_i)$ is conjugate to rotation by $2\pi a_i'/n_i = 2\pi a_i/n$ for each $i$.

By Remark \ref{rem: primeness of links} and Lemma \ref{lemma: basic top consequences of cpaf}, the link $L$ is prime and its exterior is irreducible. Hence by Theorem \ref{thm: locbc intro}, the branched cover associated to the epimorphism $\pi_1(X(L))\rightarrow \mathbb{Z}/n$ which sends $\mu_i$ to $a_i$ (mod $n$) is left-orderable. This completes the proof.
\end{proof}

\begin{remark}
\label{remark: thm 1.1 more is true}
The conclusion of Theorem \ref{thm: lo cbcs intro} also holds for epimorphism $\psi: \pi_1(X(L)) \rightarrow \mathbb{Z}/n$ as long as $n_i |\mu_i\cdot \delta_i(\Phi_0)| \geq 2$, where $n_i$ is the order of $\psi(\mu_i)$ in $\mathbb{Z}/n$. That is, if for some $i$, $|\mu_i\cdot \delta_i(\Phi_0)| \geq 2$, then the conclusion holds even if $n_i = 1$, i.e. $\psi(\mu_i) = 0$ in $\mathbb{Z}/n$. 
\end{remark}

Recall Corollary \ref{cor: $LO$ branched cover with all orientations intro}, which we restate here for the reader's convenience.

\begin{customcor}{\ref{cor: $LO$ branched cover with all orientations intro}}
Let $L = K_1 \cup \cdots \cup K_m$ be a link in an integer homology $3$-sphere $W$ whose complement admits a pseudo-Anosov flow none of whose degeneracy loci are meridional. Then $\pi_1(\Sigma_n(L^\mathfrak{o}))$ is left-orderable for all $n \geq 2$ and all orientations $\mathfrak{o}$ on $L$, 
\end{customcor}

This corollary partially confirms a rather surprising consequence of known results on the left-orderability of the fundamental groups of the $n$-fold cyclic branched covers $\Sigma_n(L)$ of prime oriented links $L$ in $S^3$. Specifically, these results are consistent with the possibility that every such link satisfies (exactly) one of the following.

\begin{itemize}
\setlength\itemsep{0.7em}
\item $\pi_1(\Sigma_n(L))$ is left-orderable for all $n \ge 2$;
\item $\pi_1(\Sigma_n(L))$ is non-left-orderable for all $n \ge 2$;
\item $\pi_1(\Sigma_n(L)) \mbox{ is } \left\{ \begin{array}{ll} 
\mbox{non-left-orderable}  & \mbox{for } 2 \le n \le N \\ 
\mbox{left-orderable} & \mbox{for } n > N  
\end{array} \right\}$  for some integer $N$, $2 \leq N \leq 5$.
\end{itemize}
A consequence of these statements would be 
\begin{equation}
\label{eqn: sigma2 implies sigman}
\mbox{if $\pi_1(\Sigma_2(L))$ is left-orderable then $\pi_1(\Sigma_n(L))$ is left-orderable for all $n \ge 2$}
\end{equation}
This seems puzzling at first sight inasmuch as $\Sigma_2(L)$ is independent of the orientation on $L$ while this is not true for $\Sigma_n(L)$, $n \ge 3$. Thus (\ref{eqn: sigma2 implies sigman}) predicts that if $\pi_1(\Sigma_2(L))$ is left-orderable then so is $\pi_1(\Sigma_n(L^\mathfrak{o}))$ for all orientations $\mathfrak{o}$ on $L$. In this direction, it follows from \cite[Theorem 1.1]{BRW05} that if $\pi_1(\Sigma_2(L))$ is left-orderable, then $\pi_1(\Sigma_{2n}(L^\mathfrak{o}))$ is left-orderable for all $n \ge 1$ and all orientations $\mathfrak{o}$ on $L$, though ostensibly it says nothing about odd order cyclic branched covers. However, Corollary \ref{cor: $LO$ branched cover with all orientations intro} gives conditions under which (\ref{eqn: sigma2 implies sigman}) does indeed hold. 

\subsection{Examples of links with non-meridional pseudo-Anosov flows}
\label{subsec; egs}
Next, we give some examples of links whose complements admit pseudo-Anosov flows with non-meridional degeneracy loci to which, therefore, Theorem \ref{thm: lo cbcs intro} applies.

\begin{thm} 
\label{thm: hyperbolic fibre sigman $LO$ intro}
Let $L = K_1 \cup \cdots \cup K_m$ be a hyperbolic link in an integer homology $3$-sphere $W$ which can be oriented to be a fibred link whose monodromy has a non-zero fractional Dehn twist coefficient on each boundary component of the fibre. Then the fundamental group of any $n$-fold cyclic branched cover of $L$, $n \geq 2$, is left-orderable. In particular, $\pi_1(\Sigma_n(L^{\mathfrak{o}}))$ is left-orderable for all $n\geq 2$ and all orientations $\mathfrak{o}$ on $L$.
\end{thm}

\begin{proof}
Since $L$ is hyperbolic, its monodromy is freely isotopic to a pseudo-Anosov homeomorphism \cite{Thurston98}. Hence, the suspension flow of the pseudo-Anosov homeomorphism gives rise to  a pseudo-Anosov flow on $C(L)$. Moreover, because the fractional Dehn twist coefficient of the monodromy on each boundary component of the fibre is nonzero, the degeneracy loci of the suspension flow are non-meridional. The theorem now follows from Theorem \ref{thm: lo cbcs intro}.
\end{proof}

The condition on the fractional Dehn twist coefficients in Theorem \ref{thm: hyperbolic fibre sigman $LO$ intro} is necessary. For example, let $L_k$ be the 2-bridge link corresponding to the continued fraction $[2,2,...,2]$ of length $k$, where $k \ge 2$. (The number of components on $L_k$ is 1 if $k$ is even and 2 if $k$ is odd; $L_2$ is the figure eight knot.) Then $L_k$ is hyperbolic and fibred, but $\Sigma_n(L_k)$ has non-left-orderable fundamental group for all $n \ge 2$ (\cite{MV2002}, \cite{IT2020}, \cite{BGW13}; see \S \ref{subsec: egs}).

The family of fibred strongly quasipositive links arises naturally as the set of bindings of open books which carry the tight contact structure on the 3-sphere (\cite{Hed10}). Topologically, Giroux's stabilization theorem characterizes the family as the set of fibred links whose fibre surface can be transformed into a plumbing of positive Hopf bands by a finite sequence of such plumbings \cite{Giroux02, Rudolph98}.
For this family, Theorem \ref{thm: hyperbolic fibre sigman $LO$ intro} implies:

\begin{cor}
\label{cor: hyperbolic fibred sqp links intro}
Suppose that $L$ is a hyperbolic link in $S^3$ which can be oriented to be fibred and strongly quasipositive. Then the fundamental group of any $n$-fold cyclic branched cover of $L$, $n \geq 2$, is left-orderable. In particular, $\pi_1(\Sigma_n(L^{\mathfrak{o}}))$ is left-orderable for all $n\geq 2$ and all orientations $\mathfrak{o}$ on $L$.  
\end{cor}

\begin{proof}
Let $\mathfrak{o}$ be an orientation on $L$ for which $L^\mathfrak{o}$ is fibred and strongly quasipositive. Since $L^\mathfrak{o}$ is strongly quasipositive, the open book associated to its fibring carries the standard tight contact structure (\cite{Hed10}). It is then shown in \cite[Theorem 1.1 and Proposition 3.1]{HKMI} that the fractional Dehn twist coefficients of its monodromy on the boundary components of its fibre are either all positive or all negative. Hence, the conclusion holds by Theorem \ref{thm: hyperbolic fibre sigman $LO$ intro}.
\end{proof}

Next we show how pseudo-Anosov flows with non-meridional degeneracy loci often exist on the complements of pseudo-Anosov closed braids.

Let $D_w$ denote the $w$-punctured $2$-disk. There is a natural identification of $\mbox{MCG}(D_w)$ with the $w$-strand braid group $B_w$ (see e.g. \cite[Chapter 9]{FM12}) and as such we can associate a fractional Dehn twist coefficient $c(b) \in \mathbb Q$ to each $b \in B_w$.

\begin{thm}
\label{thm: closed braids}
Let $b \in B_w$ be a pseudo-Anosov braid and $L=\hat{b}$. Suppose that  $c(b)$ is neither $0$ nor the reciprocal of a non-zero integer. Then the fundamental group of any $n$-fold cyclic branched cover of $L$, $n \geq 2$, is left-orderable.  In particular, $\pi_1(\Sigma_n(L^{\mathfrak{o}}))$ is left-orderable for all $n\geq 2$ and all orientations $\mathfrak{o}$ on $L$. 
\end{thm}

\begin{proof} 
By assumption $b$, thought of as a mapping class, is freely isotopic to a pseudo-Anosov homeomorphism $\beta$ of $D_w$.

The interior of the mapping torus of $\beta$ is the complement of an $m+1$ component link in $S^3$ consisting of $L$ and the braid axis $A$, while the restriction of the suspension flow $\Phi_0$ of $\beta$ to the complement of $L\sqcup A$ is pseudo-Anosov. Express the degeneracy locus $\delta_A(\Phi_0)$ of $\Phi_0$ on the boundary of a tubular neighbourhood of $A$ homologically as 
$$\delta_A(\Phi_0) = p\mu' + q\lambda'$$
where $\mu'$ and $\lambda'$ are meridional and longitudinal classes of $A$. Our hypotheses on $c(b)$ imply that $|q| > 1$ and therefore $\mu'$ intersects $\delta_A(\Phi_0)$ at least twice. Hence  there is a pseudo-Anosov flow $\bar \Phi_0$ on the complement of $L$ obtained by extending $\Phi_0$ over a tubular neighbourhood of $A$. It is clear that the degeneracy loci of $\bar \Phi_0$ are non-meridional, so the desired conclusion follows from Theorem \ref{thm: lo cbcs intro}.
\end{proof} 

\begin{remark}
Theorem \ref{thm: closed braids} extends with the same proof to the closures of pseudo-Anosov braids in open book decompositions of integer homology $3$-spheres.  
\end{remark}

It was shown in \cite[Theorem 1.9]{BH19} that given a pseudo-Anosov braid $b$ on an odd number of strands, if $|c(b)|\geq 2$, then all even order cyclic branched covers of $\hat{b}$ have left-orderable fundamental groups and admit co-oriented taut foliations. The following corollary substantially improves the left-orderability part of this result.

\begin{cor}
\label{cor: fdtc braid}
Let $L = \hat{b}$ be an hyperbolic link in $S^3$, where the  fractional Dehn twist coefficient of $b \in B_w$ satisfies $|c(b)|>1$. Then the fundamental group of any $n$-fold cyclic branched cover of $L$, $n \geq 2$, is left-orderable.  In particular, $\pi_1(\Sigma_n(L^{\mathfrak{o}}))$ is left-orderable for all $n\geq 2$ and all orientations $\mathfrak{o}$ on $L$. 
\end{cor}

\begin{proof}
Since $L$ is hyperbolic and $|c(b)|>1$, \cite[Theorem 8.4]{IK07} implies that $b$ is pseudo-Anosov braid. The corollary now follows from Theorem \ref{thm: closed braids} since the condition $|c(b)| > 1$ implies that $c(b)$ is neither zero nor the reciprocal of a non-zero integer.  
\end{proof}

\begin{remark}
In \S \ref{subsec: egs} we give infinitely many examples of braids $b$ with $c(b) = 1$ such that $L = \hat b$ is a $2$-component hyperbolic link with $\pi_1(\Sigma_n(L))$ left-orderable for all $n \geq 3$, while if $L^\mathfrak{o}$ is the oriented link obtained by reversing the orientation of one of components of $L$ then $\pi_1(\Sigma_n(L^\mathfrak{o}))$ is non-left-orderable for all $n \geq 2$. This shows that Corollary \ref{cor: fdtc braid} can fail quite dramatically if the condition $|c(b)| > 1$ is relaxed. 
\end{remark}

Since adding a positive full twist to a braid $b$ will increase its fractional Dehn twist coefficient by $1$, it is easy to obtain examples to which Corollary \ref{cor: fdtc braid} can be applied. More precisely, recall that the centre of $B_w$ is generated by the braid $C_w$ corresponding to a positive full Dehn twist along $\partial D_w$. Then,

\begin{prop}
\label{prop: example of braids}
Assume that $b$ is a pseudo-Anosov braid in $B_w$ and for $k \in \mathbb Z$, let $L_k$ be the closure of the braid $C_w^k b$. Suppose
\begin{enumerate}[leftmargin=*]
\setlength\itemsep{0.3em}
\item[{\rm (1)}] $k \neq -c(b), -c(b)\pm 1$ when $c(b)\in \mathbb{Z};$ 
\item[{\rm (2)}]  $k \neq -\lfloor c(b)\rfloor, -\lfloor c(b)\rfloor-1$ when $c(b) \notin \mathbb{Z}$.  
\end{enumerate}
Then the fundamental group of any $n$-fold cyclic branched cover of $L_k$, $n \geq 2$, is left-orderable.  In particular, $\pi_1(\Sigma_n(L_k^{\mathfrak{o}}))$ is left-orderable for all $n\geq 2$ and all orientations $\mathfrak{o}$ on $L_k$. 
\end{prop}

\begin{proof}
Since $C_w$ corresponds to a positive full Dehn twist along $\partial D_w$, $C_w^k b$ is freely isotopic to $b$ and is therefore also pseudo-Anosov. Hence $\widehat{C_w^k b}$ is hyperbolic when $|c(C_w^k b)| > 1$ (\cite[Theorem 8.4]{IK07}). We have $c(C_w^k b) = k+ c(b)$ by the definition of the fractional Dehn twist coefficient of a braid (cf. \cite[Proposition 2.7]{KR13}), and the reader will verify that the condition $|c(C_w^k b)| > 1$ corresponds to the conditions in the statements of (1) and (2) of the proposition. An application of Corollary \ref{cor: fdtc braid} completes the proof.
\end{proof}

Recall that a braid $b$ in $B_w$ is {\it quasipositive} if it is represented by a braid word of the form 
\begin{displaymath}
 b = \prod_{i=1}^s w_i \sigma_{k_i} w_i^{-1}
\end{displaymath} 
where $\sigma_1, \cdots, \sigma_{w-1}$ are the standard generators of $B_w$ (\cite{Rudolph83}). 

\begin{lemma}
\label{lem: fdtc quasi positive braids}
Let $b$ be a pseudo-Anosov quasipositive braid in $B_w$. Then $c(b)>0$.  
\end{lemma} 
\begin{proof}
Let $S$ be the double branched cover of $D_w$ over $w$ points contained in $\mbox{int}(D_w)$ and $\varphi \in {\rm MCG(S)}$ be the lift of $b$. Then $\varphi$ is pseudo-Anosov and $c(\varphi) = \frac{c(b)}{2}$ (cf. \cite[Lemma 10.1,  Lemma 10.2]{BH19}). It is known that each $\sigma_i$ lifts to a positive Dehn twist along a simple closed curve in $S$ \cite[\S 9.4]{FM12}, and the same is true of any conjugate of $\sigma_i$. Then the quasipositivity of $b$ implies that $\varphi$ is a product of positive Dehn twists. Thus $\varphi$ is right-veering and hence $c(b) = 2c(\varphi) > 0$ (\cite[Proposition 3.1]{HKMI}).
\end{proof}

Theorem \ref{thm: example of braids} below follows immediately from Proposition \ref{prop: example of braids} and Lemma \ref{lem: fdtc quasi positive braids}.

\begin{thm}
\label{thm: example of braids}
Let $b$ be a pseudo-Anosov quasipositive braid in $B_w$ and $L_k$ the braid closure of $C_w^k b$, $k \in \mathbb{Z}$. Then for any $k \geq 1$, the fundamental group of any $n$-fold cyclic branched cover of $L_k$, $n \geq 2$, is left-orderable.  In particular, given any $k\geq 1$, $\pi_1(\Sigma_n(L_k^{\mathfrak{o}}))$ is left-orderable for all $n\geq 2$ and all orientations $\mathfrak{o}$ on $L_k$.
\end{thm}

\subsection{Dependence of the left-orderability of the fundamental groups of cyclic branched covers on link orientation}
\label{subsec: egs}
Here we show that in contrast to expectation (\ref{eqn: sigma2 implies sigman}), if $\pi_1(\Sigma_2(L))$ is not left-orderable then in general the left-orderability of $\pi_1(\Sigma_n(L))$ for $n > 2$ depends on the orientation on $L$.

\begin{figure}[ht]
\centering
\includegraphics[scale=0.6]{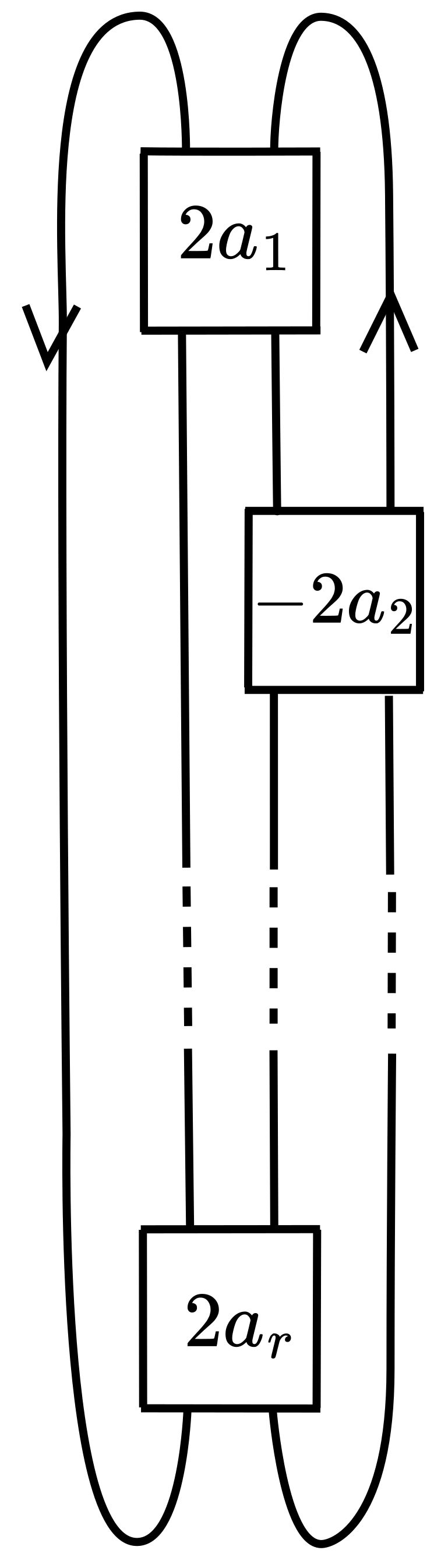}
\caption[short]{The canonical orientation of 2-bridge links}
\label{fig: canonical orientation}
\end{figure}

Given non-zero integers $a_1,a_2,...,a_r$, $r \ge 1$, let $L = L(2a_1,2a_2,...,2a_r)$ be the 2-bridge knot or link corresponding to the rational number with continued fraction $[2a_1,2a_2,...,2a_r]$. This is a knot if $r$ is even and a 2-component link if $r$ is odd. (We use the ``+'' convention for continued fractions, so for example $L(2,2)$ is the figure eight knot, corresponding to the rational number 5/2.) When $L$ has two components we give it the {\it canonical} orientation, illustrated in Figure \ref{fig: canonical orientation}. We then use $L^\mathfrak{o}$ to denote the link obtained by reversing the orientation of one of the components of $L$.

If all the $a_i$ are positive then $\Sigma_n(L) \cong \Sigma_2(L_n)$ for some alternating link $L_n$ \cite{MV2002}, \cite{IT2020}, and hence $\pi_1(\Sigma_n(L))$ is non-left-orderable for all $n \ge 2$ by \cite{BGW13}. In particular, for the $(2,2k)$-torus links $L(2k)$, $k \ge 1$, $\pi_1(\Sigma_n(L(2k))$ is non-left-orderable for all $n \ge 2$. (This was first proved by Dabkowski, Przytycki and Togha in \cite{DPT05} by a direct algebraic argument.) On the other hand, with the other orientation, $\pi_1(\Sigma_n(L(2k)^\mathfrak{o}))$ is left-orderable for all $n \ge 3$ if $k \ge 3$, and for all $n \ge 4$ if $k = 2$ \cite{BGH3}. (We remark that the orientation on $L(2k)^\mathfrak{o}$ is the one coming from its realization as the closure of the 2-braid $\sigma_1^{2k}$.) This shows that in general, if $\pi_1(\Sigma_2(L))$ is not left-orderable then the left-orderability of $\pi_1(\Sigma_n(L))$ for $n > 2$ depends on the orientation on $L$. Note that torus links are Seifert links, i.e. their exteriors are Seifert fibered. The following theorem gives hyperbolic examples.     

\begin{thm} 
\label{thm: hyp egs}
Let $L = L(2k_1,2l_1,...,2k_r,2l_r,2k_{r+1})$, where $r$, $k_i$, and $l_i$ are positive. Then, 
\begin{enumerate}[leftmargin=*]
\setlength\itemsep{0.3em}
    \item[{\rm (1)}] $\pi_1(\Sigma_n(L))$ is non-left-orderable for all $n \ge 2$.
    \item[{\rm (2)}] If $k_i = k \geq 3$ for all $i$, then $\pi_1(\Sigma_n(L^\mathfrak{o}))$ is left-orderable for all $n \geq 3$. 
    \item[{\rm (3)}] If $k_i = 2$ for all $i$, then $\pi_1(\Sigma_n(L^\mathfrak{o}))$ is left-orderable for all $n \geq 4$. 
    
    \item[{\rm (4)}] If $k_i = l_i = 1$ for all $i$ and $r$ is odd, then $\pi_1(\Sigma_n(L^\mathfrak{o}))$ is left-orderable for all $n \geq 7$. 
\end{enumerate}

\end{thm}
Note that for the links in part (4) of the theorem $L$ is fibred. 

\begin{proof}
Part (1) is discussed above.

For parts (2), (3), and (4) we use a result of Ohtsuki, Riley and Sakuma \cite{ORS08}, who constructed epimorphisms between $2$-bridge link groups which preserve oriented meridional classes with respect to the canonical orientations, and hence also the non-canonical orientations. Any such epimorphism induces epimorphisms between the fundamental groups of the associated cyclic branched covers so if the target is left-orderable, so is the domain. 

For (2) and (3) we start with $L(2k)$, $k \ge 2$. It follows from \cite[Proposition 5.1]{ORS08} that for $L$ of the form $L(2k,2l_1,2k,2l_2,...,2l_r,2k)$ there is an epimorphism $\pi_1(S^3 \setminus L) \to \pi_1(S^3 \setminus L(2k))$ as above. Then (2) and (3) follow from the remarks in the paragraph immediately preceding Theorem \ref{thm: hyp egs}.

For (4) we start with $L(2,2,2)$. Lemma \ref{lemma: lon7} below shows that  $\pi_1(\Sigma_n(L(2,2,2)^{\mathfrak{o}}))$ is left orderable for all $n \ge 7$. It follows from \cite[Proposition 5.1]{ORS08} that if 
$$L = L(2,2,2,2l_1,2,2,2,2l_2,...,2,2,2,2l_r,2,2,2)$$ then  $\pi_1(\Sigma_n(L^{\mathfrak{o}}))$ is left-orderable for $n \geq 7$. In particular, taking all $l_i = 1$ we get $L(2,2,2,...,2)$ of any length congruent to 3 (mod 4).
\end{proof}

\begin{figure}[ht]
    \includegraphics[scale=0.5]{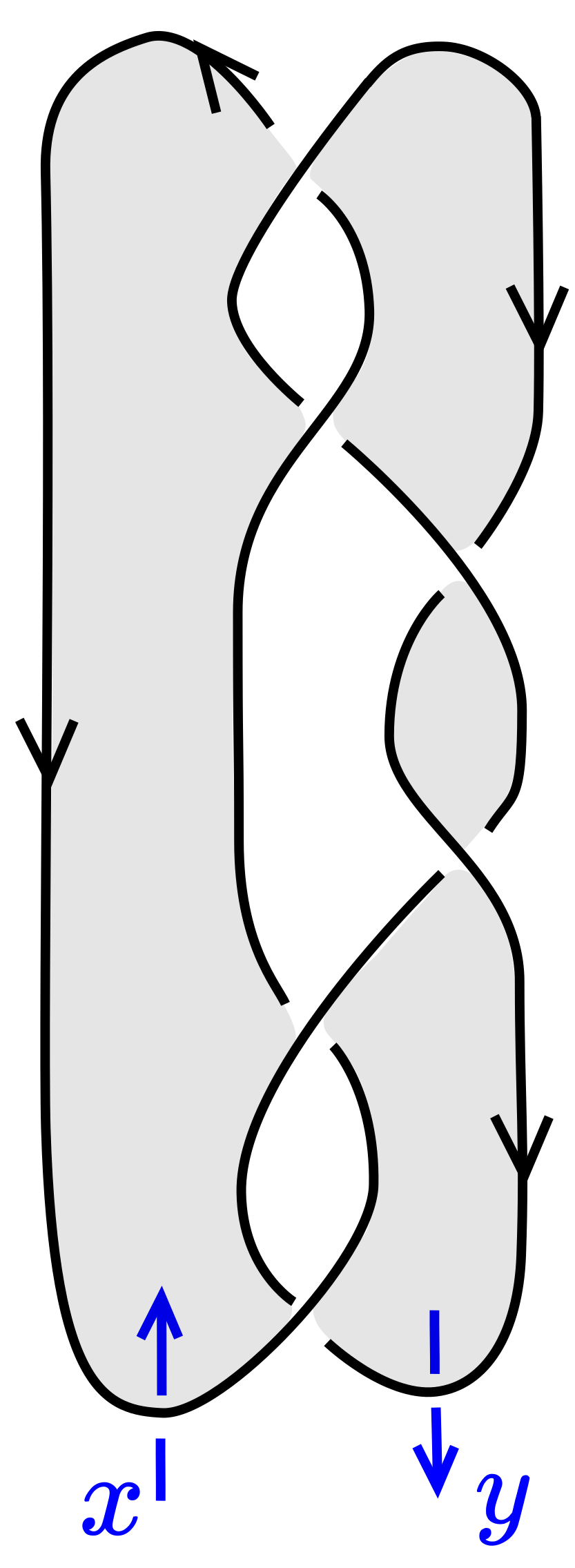}
    \caption{The 2-bridge link $L(2, 2, 2)^{\mathfrak{o}}$ with the non-fibered orientation}
    \label{fig: two bridge}
\end{figure}

\begin{lemma}
\label{lemma: lon7}
$\pi_1(\Sigma_n(L(2,2,2)^\mathfrak{o}))$ is left-orderable for all $n \geq 7$. 
\end{lemma}

\begin{proof}
 We write down a presentation of the link group using the link diagram in Figure \ref{fig: two bridge}:  
    \begin{displaymath}
        \pi_1(X(L)) = <x, y : w x = x w>,
    \end{displaymath}
    where $x$ and $y$ are meridional generators in the standard Wirtinger presentation shown in Figure \ref{fig: two bridge} and $w = yxy^{-1}x^{-1}yxyx^{-1}y^{-1}xy$. 

    Define $\rho_\theta: \pi_1(X(L))\rightarrow SL(2,\mathbb{C})$ by setting
    \begin{displaymath}
        \rho(x) = \begin{pmatrix}
            m & 1 \\
            0&  m^{-1}
        \end{pmatrix}, \quad  \rho(y) = \begin{pmatrix}
            m & 0 \\
            s &  m^{-1}
        \end{pmatrix}
    \end{displaymath}
    where $m = e^{i\theta}$ and $s = 3 - 4\cos^2(\theta)$, $\theta\neq 0$. One can easily verify that $$\rho_\theta(w) = \begin{pmatrix}
    - \cos(3\theta) + i \sin(3\theta) & 1 + 2 \cos(2\theta) \\
    0 & -\cos(3\theta) - i \sin(3\theta)
\end{pmatrix}$$
and hence $\rho_\theta(w x) = \rho_\theta(x w)$. So $\rho_\theta$ defines an $SL(2,\mathbb{C})$-representation of $\pi_1(X(L))$ for each $\theta\neq 0$. Also note that $\rho_\theta(y^{-1}xy) = \rho_\theta(xyx^{-1})$. 

By \cite[p. 786]{Kho} (also see \cite[Theorem 4.3]{Hu15}), when $s(\theta)<0$ or $s(\theta) > 4\sin^2(\theta)$, the representation $\rho_\theta$ is conjugate to an $SL(2,\mathbb{R})$-representation, denoted by $\rho'_\theta$, for which $\rho'_\theta(x)$ is conjugate to the rotation by $\theta$. 

Hence, letting $\theta = \frac{\pi}{n}$, for each $n > 6$, we have an $SL(2,\mathbb{R})$-representation $\rho'_\theta$, where $\rho'(x)$ is conjugate to rotation by $\frac{\pi}{n}$. Since $\rho_\theta(y^{-1}xy) = \rho_\theta(xyx^{-1})$, we have $\rho'_\theta(y^{-1}xy) = \rho'_\theta(xyx^{-1})$, and therefore, $\rho'_\theta(y)$ is also conjugate to a rotation by $\frac{\pi}{n}$. 

By projecting the representation to a $PSL(2, \mathbb{R})$-representation, for each $n>6$ we have a $PSL(2,\mathbb{R})$ representation $\rho$ of $\pi_1(X(L))$ such that $\rho(x)$ and $\rho(y)$ are both conjugate to the rotation by $\frac{2\pi}{n}$. It follows from Theorem \ref{thm: locbc intro} that $\pi_1(\Sigma_n(L))$ is left-orderable for all $n > 6$.
\end{proof}

In Section 7 we prove an analog of Theorem \ref{thm: hyp egs} where the property ``has left-orderable fundamental group'' is replaced by ``is not an $L$-space'', consistent with the $L$-space Conjecture.

We now show that there are hyperbolic 2-bridge links which show that Corollary \ref{cor: fdtc braid} is best possible. We first note the following, which will also be used in Section 7.

\begin{lemma}
$L^\mathfrak{o}$ is strongly quasipositive.
\end{lemma}

\begin{proof}
It is clear that $L^\mathfrak{o}$ is special alternating (i.e. one of the chessboard surfaces in an alternating diagram is orientable; cf. Figure \ref{fig: two bridge}), therefore positive, and hence strongly quasipositive by \cite{Rudolph99}.  
\end{proof}

\begin{figure}[ht]
\includegraphics[scale=0.5]{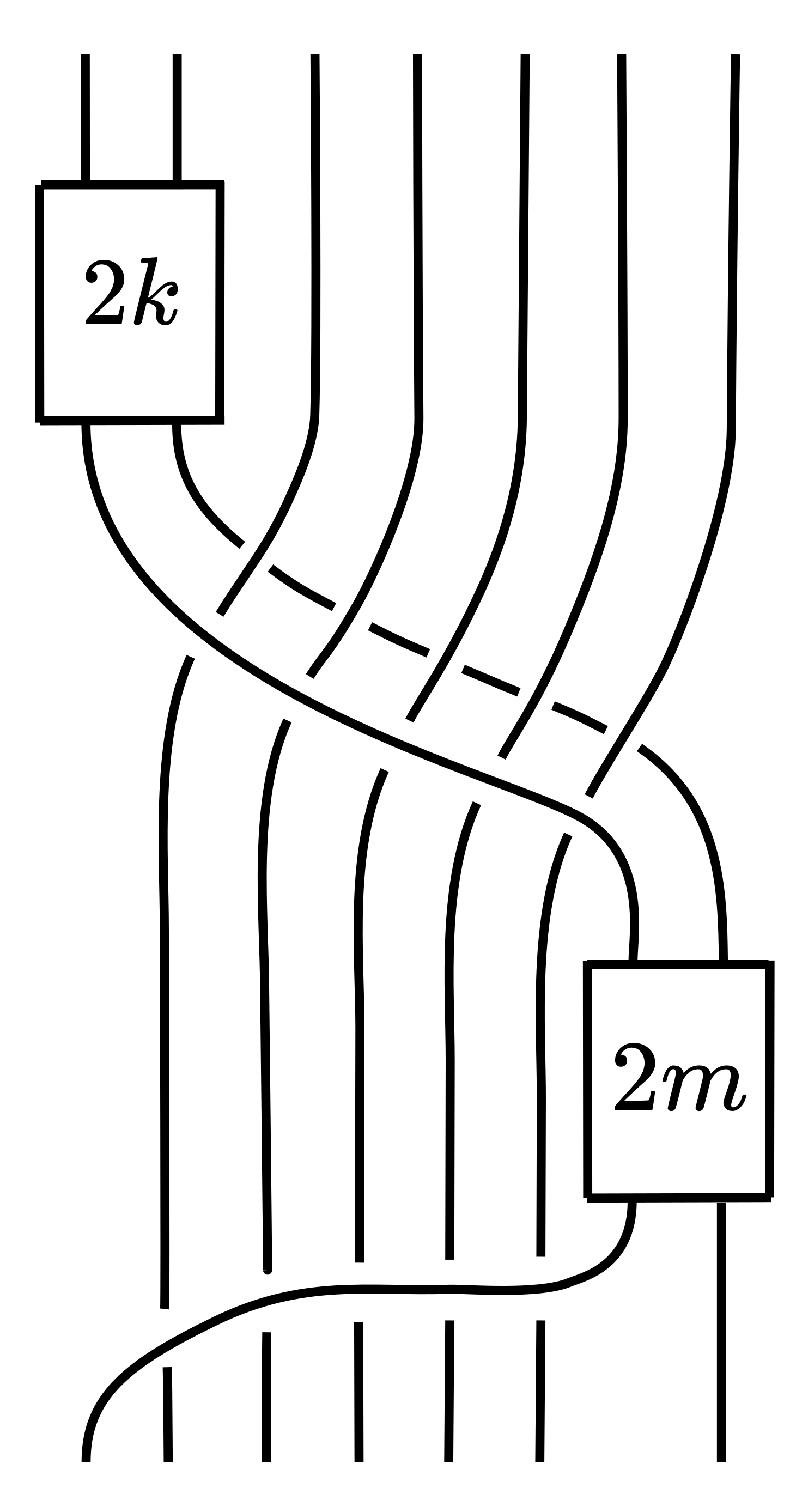}
\caption[]{A closed braid presentation of $L(2k, 2l, 2m)^{\mathfrak{o}}$ with $l=5$.}
\label{fig: L as a closed braid}
\end{figure}

 
 In the case $r = 1$ we identify an explicit strongly quasipositive braid whose closure is $L^\mathfrak{o}$. Write $L = L(2k,2l,2m), \, k,l,m \ge 1$. Then it is straightforward to verify that $L^\mathfrak{o} = \hat{b}$ where $b$ is the strongly quasipositive index $(l+2)$ braid
$$\sigma_1^{2k}(\sigma_2...\sigma_{l+1})(\sigma_1^{-1}...\sigma_l^{-1})\sigma_{l+1}^{2m}(\sigma_l...\sigma_1) = \sigma_1^{2k}(\sigma_2...\sigma_{l+1})a_{1,l+2}^{2m}$$
shown in Figure \ref{fig: L as a closed braid}, which illustrates the case $l = 5$. Here $a_{1,l+2} = (\sigma_1^{-1}\sigma_2^{-1}...\sigma_l^{-1})\sigma_{l+1}(\sigma_l...\sigma_2\sigma_1)$. 
  
The following shows that the condition $c(b) > 1$ in Corollary \ref{cor: fdtc braid} cannot be relaxed.

\begin{thm}
  Let $L$ be the closure of the 3-braid $b = \sigma_1^{2k} \sigma_2 a_{13}^{2k}$, where $k \ge 3$, and let $L^\mathfrak{o}$ be obtained by reversing the orientation of one of the components of $L$. Then $c(b) = 1$, $\pi_1(\Sigma_n(L))$ is left-orderable for all $n \ge 3$,  and $\pi_1(\Sigma_n(L^\mathfrak{o}))$ is non-left-orderable for all $n \ge 2$.
\end{thm}

\begin{proof}
  Taking $l = 1$ in the $(l+2)$-braid $b$ defined above we get the the 3-braid $b = \sigma_1^{2k}\sigma_2a_{13}^{2m}$. More generally, let $b(p,q,r)$ be the strongly quasipositive 3-braid $\sigma_1^p\sigma_2^qa_{13}^r, \, p,q,r \ge 1$. These braids are considered in \cite[Corollary 3.7]{BBG19-2}, which lists the strongly quasipositive 3-braids whose closures are prime, non-split, non-trivial, definite links. Among these are the braids $b(p,q,r)$, whose closures are not fibred. (The closures of the others are precisely the ADE links, i.e. the fibered strongly quasipositive links whose fiber is a plumbing of positive Hopf bands according to the tree associated to a Dynkin diagram of type A, D, or E. See e.g. \cite{BGH3}).  Taking account of the different braid word conventions used in \cite{BBG19-2} and in the present paper, \cite[Lemma 3.2(3)]{BBG19-2} shows that $b(p,q,r)$ is conjugate to $C_{3}\sigma_1^{p-1}\sigma_2^{-1}\sigma_1^{q-1}\sigma_2^{-1}\sigma_1^{r-1}\sigma_2^{-1}$. In \cite[Proof of Theorem 1.11]{BH19} it is shown that the fractional Dehn twist coefficient $c(b(p,q,r)) = 1$. In particular, taking $p = 2k, q = 1, r = 2m$, the braid $b = \sigma_1^{2k}\sigma_2a_{13}^{2m}$ has $c(b) = 1$. Reorienting $L^\mathfrak{o} = \hat{b}$ gives $L = L(2k,2,2m)$, and $\Sigma_n(L)$ has a non-left-orderable fundamental group for all $n \ge 2$ by Theorem \ref{thm: hyp egs}(1). On the other hand, taking $k = m \ge 3$, $\pi_1(\Sigma_n(L^\mathfrak{o}))$ is left-orderable for all $n \ge 3$ by Theorem \ref{thm: hyp egs}(2).
\end{proof}
  
\subsection{Application to degeneracy loci}
\label{subsec: meridional degeneracy loci}
As mentioned in the introduction, Gabai and Mosher have independently shown that pseudo-Anosov flows exist on the complement of any hyperbolic link in a closed, orientable 3-manifold. More precisely, they show that given a finite depth taut foliation $\mathcal{F}$ on a compact, connected, orientable, hyperbolic $3$-manifold $M$ with non-empty boundary consisting of tori, there is a pseudo-Anosov flow on the interior of $M$ which is almost transverse to $\mathcal{F}$ (cf. \cite[Theorem C(3)]{Mosher96}). Unfortunately, no proof has been published, though Landry and Tsang have recently produced the first of several planned articles which will provide a demonstration. See \cite{LT}.
  
However, the degeneracy loci of these flows are difficult to determine in general. Corollary \ref{cor: $LO$ branched cover with all orientations intro} implies that if $L$ is a link in an integer homology $3$-sphere which admits an orientation $\mathfrak{o}$ for which $\Sigma_n(L^\mathfrak{o})$ has a non-left-orderable fundamental group for some $n\geq 2$, then the degeneracy loci of a pseudo-Anosov flow on the link's complement cannot all be non-meridional. Specializing to knots we have: 
  
  \begin{cor}
  \label{cor: degeneracy loci $LO$ branched covers intro}
  Let $K$ be a knot in an integer homology $3$-sphere. If $\Sigma_n(K)$ has a non-left-orderable fundamental group for some $n\geq 2$, then the degeneracy locus of any pseudo-Anosov flow on the complement of $K$ is meridional.
  \qed
  \end{cor}
  For instance, the $2$-fold cyclic branched covers of alternating knots have non-left-orderable fundamental groups (\cite[Theorem 4]{BGW13}), so we obtain:   
  
  \begin{cor}
  \label{cor: degeneracy loci hyp alt knots}
  The degeneracy locus of any pseudo-Anosov flow on the complement of an alternating knot is meridional.  
  \qed
  \end{cor}

\section{Connections with the \texorpdfstring{$L$}{L}-space Conjecture}
\label{sec: applications}
  
One of the motivations for studying the left-orderability of $3$-manifold groups is the $L$-space Conjecture \cite{BGW13}, \cite{Juhasz2015}, which asserts that for a prime, closed, orientable $3$-manifold $M$ the following are equivalent:
\begin{enumerate}[leftmargin=*]
\setlength\itemsep{0.3em}
\item $M$ is not an $L$-space in the sense of Heegaard Floer homology;
\item $\pi_1(M)$ is left-orderable;
\item $M$ admits a co-orientable taut foliation.
\end{enumerate}
It is known that (3) implies (1) (\cite{OS04}, \cite{KR19}, \cite{Bow16}); all other implications are open in general. 

In this context it is interesting to compare Corollary \ref{cor: hyperbolic fibred sqp links intro}, which we restate below, with \cite[Theorem 1.1]{BBG19-1}.

\begin{customcor}{\ref{cor: hyperbolic fibred sqp links intro}}
Suppose that $L$ is a hyperbolic link in $S^3$ which can be oriented to be fibred and strongly quasipositive. Then the fundamental group of any $n$-fold cyclic branched cover $\Sigma_{\psi}(L)$ of $L$, $n \geq 2$, is left-orderable. In particular, $\pi_1(\Sigma_n(L^{\mathfrak{o}}))$ is left-orderable for all $n\geq 2$ and all orientations $\mathfrak{o}$ on $L$.  
\end{customcor}

The following is part of \cite[Theorem 1.1]{BBG19-1}.
  
\begin{thm}
\label{thm: bbg1}
If $L$ is fibred and strongly quasipositive then $\Sigma_n(L)$ is not an $L$-space for all $n \geq 6$.
\end{thm}
In this generality, Theorem \ref{thm: bbg1} is best possible: for the trefoil $T(2,3)$, $\Sigma_n(T(2,3))$ has finite fundamental group, and therefore is an $L$-space, for $2 \leq n \leq 5$. However, an obvious problem in light of Corollary \ref{cor: hyperbolic fibred sqp links intro} is to show that if $L$ is assumed to be hyperbolic then the conclusion of Theorem \ref{thm: bbg1} holds for all $n \geq 2$. 

Another obvious difference between Corollary \ref{cor: hyperbolic fibred sqp links intro} and Theorem \ref{thm: bbg1} is the independence in Corollary \ref{cor: hyperbolic fibred sqp links intro} of the particular $n$-fold cyclic branched cover $\Sigma_\psi(L)$. A challenge is to show that this holds with ``$\pi_1(\Sigma_\psi(L))$ is left-orderable" replaced by ``$\Sigma_\psi(L)$ is not an $L$-space".

We next discuss the three properties in the $L$-space Conjecture for the cyclic branched covers of $L$-space knots. These are prime \cite{Krc15}, fibred \cite{Ni07}, and strongly quasipositive \cite{Hed10}.

First, the situation is completely understood for torus knots \cite{GL14} (see also \cite{BGH3}): $\Sigma_n(K)$ has left-orderable fundamental group, is not an $L$-space, and admits a co-orientable taut foliation if and only if $\pi_1(\Sigma_n(K))$ is infinite. In particular, $\Sigma_n(K)$ has all three properties for all $n \geq 2$ if and only if $K$ is not $T(3, 4), T(3, 5)$, or $T(2, 2q+1)$ for some $q \geq 1$.

Second, if $K$ is any satellite knot, then for all $n \geq 2$ $\Sigma_n(K)$ has left-orderable fundamental group, is not an $L$-space, and, if the companion is fibred, admits a co-orientable taut foliation \cite{BGH21}. In particular this applies to $L$-space knots. 

This leaves the hyperbolic case, where we have the following.

\begin{thm}
\label{thm: hyp case}
Let $K$ be a hyperbolic $L$-space knot. Then,
\begin{enumerate}[leftmargin=*]
\setlength\itemsep{0.3em}
    \item[{\rm (1)}] $\pi_1(\Sigma_n(K))$ is left-orderable for all $n \geq 2$.
    \item[{\rm (2)}] $\Sigma_n(K)$ is not an $L$-space for all $n \geq 3$.
    \item[{\rm (3)}] $\Sigma_n(K)$ admits a co-orientable taut foliation for all $n \geq 4g(K) - 2$.
\end{enumerate}

\end{thm}

\begin{proof}
Part (1) follows from Corollary \ref{cor: hyperbolic fibred sqp links intro}. Part (2) follows from \cite[Corollary 1.4]{BBG19-1} and \cite{FRW22}. The former says that the conclusion holds for $n \geq 4$ , and for $n = 3$ unless $g(K) = 2$, while the latter says that the only $L$-space knot with genus $2$ is $T(2,5)$. Part (3) follows from \cite{BH19}. 
\end{proof}

Combining part (1) of Theorem \ref{thm: hyp case} with the results for torus knots and satellite knots discussed above gives the following.

\begin{cor}
\label{cor: l-space knots branched covers}
If $K$ is an $L$-space knot then $\pi_1(\Sigma_n(K))$ is left-orderable for all $n \ge 2$ if and only if $K$ is not $T(3, 4), T(3, 5)$, or $T(2, 2q+1)$ for some $q \geq 1$.
\end{cor}

Similarly, using part (2) of Theorem \ref{thm: hyp case}, we get the analogous statement with ``$\pi_1(\Sigma_n(K))$ is left-orderable" replaced by ``$\Sigma_n(K)$ is not an $L$-space", provided $n \ge 3$. This leaves open the interesting question, due to Allison Moore, asking whether the double branched cover of a hyperbolic $L$-space knot can ever be an $L$-space. 

The discussion at the beginning of \S \ref{subsec: lo and link orientation} of the $n$-fold cyclic branched covers $\Sigma_n(L)$ of prime oriented links $L$ in $S^3$ applies with the property ``has left-orderable fundamental group" replaced by ``is not an $L$-space". The following theorem is a version of Theorem \ref{thm: hyp egs} for the latter property, describing how the analog of (\ref{eqn: sigma2 implies sigman}) can fail if $\Sigma_2(L)$ is an $L$-space.

\begin{thm} 
\label{thm: hyp egs nls}
Let $L$ and $L^\mathfrak{o}$ be the $2$-bridge links in Theorem \ref{thm: hyp egs}. Then, 
\begin{enumerate}[leftmargin=*]
\setlength\itemsep{0.3em}
\item[{\rm (1)}] $\Sigma_n(L)$ is an $L$-space for all $n \geq 2$.  
\item[{\rm (2)}] If some $k_i \geq 3$, then $\Sigma_n(L^\mathfrak{o})$ is not an $L$-space for all $n \geq 3$. 
\item[{\rm (3)}] If some $k_i = 2$, then $\Sigma_n(L^\mathfrak{o})$ is not an $L$-space for all $n \geq 4$. 
\item[{\rm (4)}]  If $k_i = 1$ for all $i$, then $\Sigma_n(L^\mathfrak{o})$ is not an $L$-space for all $n \geq 2\pi/\arccos(l/l+1)$, where $l = \min\{l_i \; | \; 1 \leq i \leq r\}$. Hence if some $l_i = 1$ then $\Sigma_n(L^\mathfrak{o})$ is not an $L$-space for all $n \geq 6$.  
\end{enumerate}
\end{thm}

Note that (4) shows that for the fibered links $L = L(2,2,...,2)$, $\Sigma_n(L^\mathfrak{o})$ is not an $L$-space for all $n \ge 6$; cf. part (4) of Theorem \ref{thm: hyp egs}.

\begin{proof}
As noted in the proof of Theorem \ref{thm: hyp egs}, $\Sigma_n(L) \cong \Sigma_2(L_n)$ where $L_n$ is alternating. It follows that $\Sigma_n(L)$ is an $L$-space \cite{OS05} for all $n \geq 2$. This proves (1).

To prove (2), (3), and (4) let $k = \max\{k_i \; | \; 1 \leq i \leq r+1\}$. 

Let $F$ be the Seifert surface for $L^\mathfrak{o}$ obtained from Seifert's algorithm applied to the link diagram shown in Figure \ref{fig: canonical orientation} but with the non-canonical orientation (cf. Figure \ref{fig: two bridge}), and let $\mathcal{S}_F$ be the associated Seifert form. For $\xi \in S^1$, let $\mathcal{S}_F(\xi)$ be the Hermitian form $(1-\xi)\mathcal{S}_F + (1-\bar{\xi})\mathcal{S}^T_F$. Choose $i$ such that $k_i = k$. Then $F$ has a subsurface $F(k)$ as shown in Figure \ref{fig: fig2}, with associated Seifert form $\mathcal{S}_{F(k)}$ and Hermitian form $\mathcal{S}_{F(k)}(\xi)$. Note that the boundary of $F(k)$ is the torus link $T(2,2k)$, with the fibered orientation.

\begin{figure}[ht]
    \includegraphics[scale=0.55]{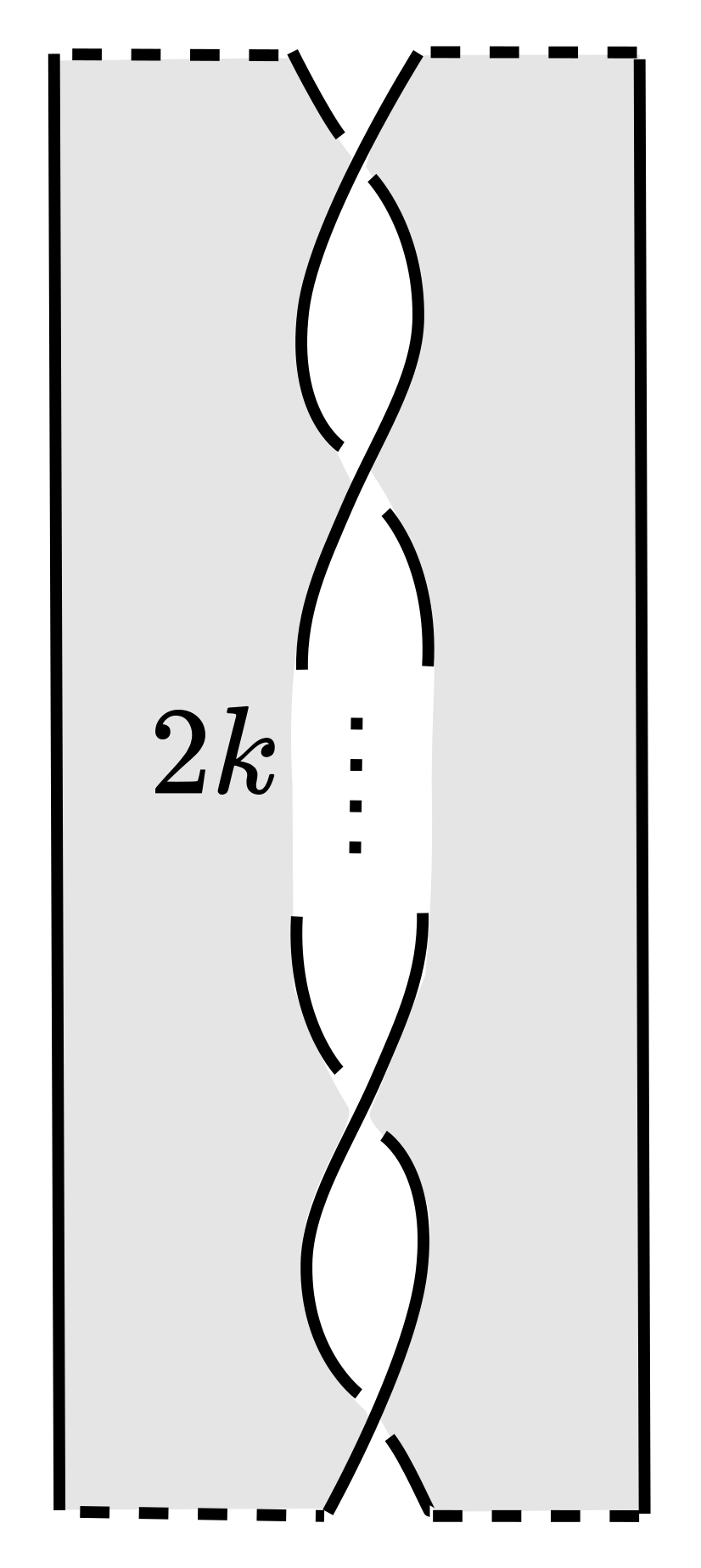}
    \caption[short]{Subspace $F(k)$ of the Seifert surface $F$ of $L^\mathfrak{o}$. There are $2k$ crossings in this diagram. }
    \label{fig: fig2}
\end{figure}


Suppose that $\Sigma_n(L^\mathfrak{o})$ is an $L$-space for some $n \ge 3$. Then by \cite[Theorem 1.1]{BBG19-1} $\mathcal{S}_F(\xi)$ is definite for $\xi \in \bar{I}_{-}(\xi_n)$, where $\xi_n = \exp(2\pi i/n)$. Since the inclusion of $F(k)$ into $F$ induces an injection of $H_1(F(k))$ into $H_1(F)$ as a direct summand, the same definiteness holds for $\mathcal{S}_{F(k)}(\xi)$. Therefore all the roots of $\Delta_{T(2,2k)}(t) = (t^{2k} - 1)/(t + 1)$ lie in $I_{+}(\xi_n)$. In particular exp$(2\pi i(k-1)/2k) \in I_{+}(\xi_n)$. Hence $(k-1)/2k < 1/n$, giving $k < n/(n-2)$. If $k \ge 3$, this implies $n < 3$, and if $k = 2$, we get $n < 4$. This proves (2) and (3).

\begin{figure}[ht]
    \centering
    \includegraphics[scale=0.5]{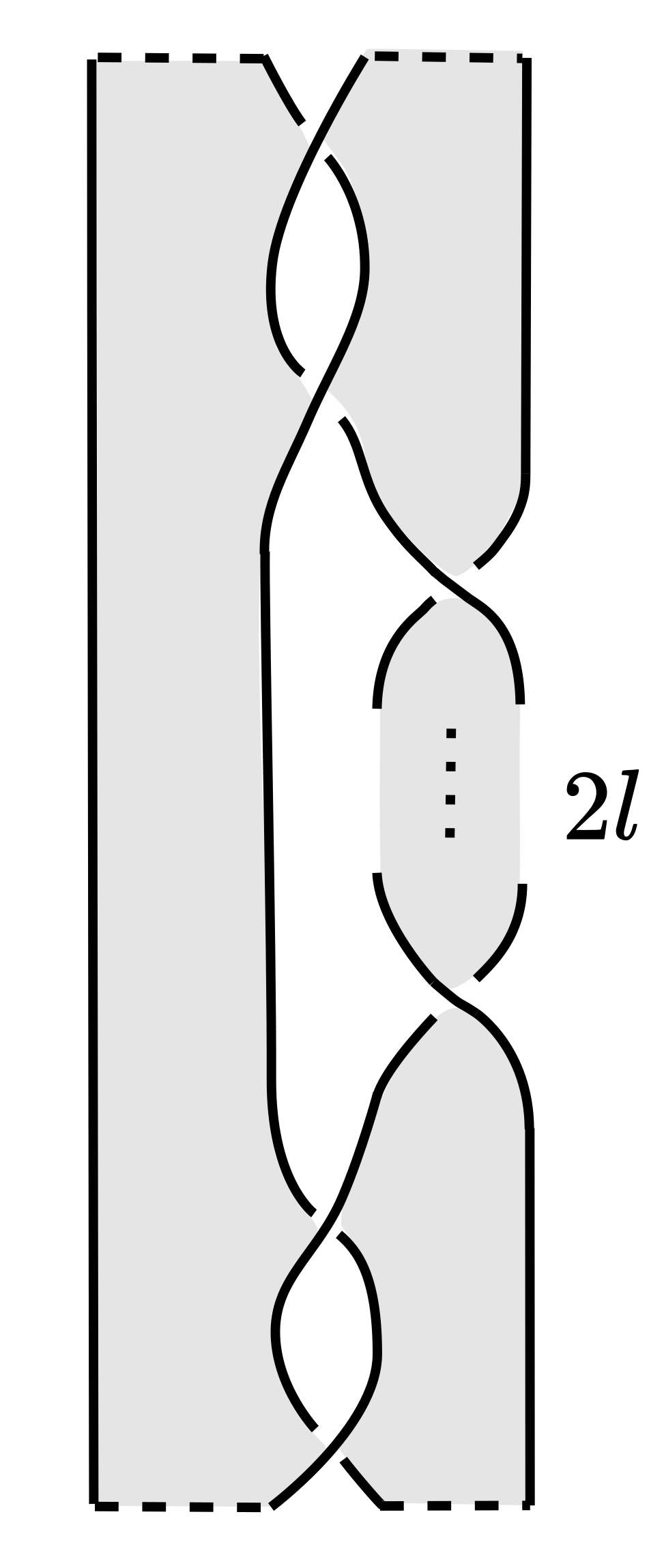}
    \caption[short]{Subsurface $F'(l)$ of the Seifert surface $F$ of $L^\mathfrak{o}$} 
    \label{fig: fig3}
\end{figure}
If $k = 1$ then $T(2,2k)$ is the Hopf link, with Alexander polynomial $(t-1)$, so instead, we choose $i$ such that $l_i = l$, and note that $F$ contains the subsurface $F'(l)$ shown in Figure \ref{fig: fig3}. As before, the Hermitian form $\mathcal{S}_{F'(l)}(\xi)$ is definite for $\xi \in \bar{I}_{-}(\xi_n)$. The link $L'(l) = \partial F'(l)$ is the 2-bridge link $L(2,2l,2)$ with the non-canonical orientation. Computing a Seifert matrix of the genus 1 surface $F'(l)$ gives $\Delta_{L'(l)}(t) = (l+1)(t-1)(1 - (2l/(l+1))t + t^2)$. The quadratic term has roots $\alpha$ and $\bar \alpha$ where $\alpha =$ exp$(i\theta)$ and 2cos$(\theta) = \alpha + \bar \alpha = 2l/(l+1).$ Hence $\theta = \,$ arccos$(l/(l+1))$. Since $\alpha \in I_{+}(\xi_n)$ we must have arccos$(l/(l+1)) < 2\pi/n$, and therefore $n < 2\pi/$arccos$(l/(l+1))$.  This proves (4).
\end{proof}

{
\footnotesize
\bibliographystyle{abbrvnat}
\bibliography{bgh_ps}
}

\vspace*{20pt}

\end{document}